\documentclass[12pt]{article}%
\usepackage{amsmath}
\usepackage{amsfonts}
\usepackage{amssymb}%
\usepackage{graphicx}
\providecommand{\U}[1]{\protect\rule{.1in}{.1in}}
\begin{document}

\title{Several Applications of Divergence Criteria in Continuous Families\\ }
\author{\ \ \ \ \ \ \ \ \ \ \ \ \ \ \ \ \ \ \ \ \ \ \ {\LARGE \ Michel Broniatowski}%
$^{a},${\LARGE Igor Vajda}$^{b,\ast}$\\$^{a}$LSTA, Universit\'{e} Paris 6 michel.broniatowski@upmc.fr\\$^{b}$Academy of Science of the Czech Republic vajda@utia.cz
\and $^{\ast}$corresponding author}
\maketitle

\begin{abstract}
This paper deals with four types of point estimators based on minimization of
information-theoretic divergences between hypothetical and empirical
distributions. These were introduced $\medskip$\newline\textbf{(i)} by Liese
\& Vajda (2006) and independently Broniatowski \& Keziou (2006), called here
\textsl{power superdivergence estimators, }\newline\textbf{(ii)}%
\textsl{\ }by\textsl{\ }Broniatowski \& Keziou (2009), called here
\textsl{power subdivergence estimators, }\newline\textbf{(iii)}\textsl{\ }%
by\textsl{\ }Basu et al. (1998), called here \textsl{power pseudodistance
estimators, }and \newline\textbf{(iv)}\textsl{\ }by Vajda (2008) called here
\textsl{R\'{e}nyi pseudodistance estimators.} $\medskip$\newline The paper
studies and compares general properties of these estimators such as
consistency and influence curves, and illustrates these properties by detailed
analysis of the applications to the estimation of normal location and
scale.\bigskip

\end{abstract}
\tableofcontents

%macros:

%%%%%%%%%%%%%%%%%%%%%%%%%%%%%%%%%%%%%%%%%%%%%%%%%

\begin{flushright}
{\tiny \today
}
\end{flushright}

%%%%%%%%%%%%%%%%%%%%%%%%%%%%%%%%%%%%%%%%%%%%%%%%%

\section{BASIC\ CONCEPTS\ AND\ RESULTS}

Let $\phi$ $:(0,\infty)\mapsto\mathbb{R}$ be twice differentiable strictly
convex function with $\phi(1)=0$ and (possibly infinite) continuous extension
to $t=0+$ denoted by $\phi(0),$ and let $\mbox{\boldmath$\Phi$}$ be the class
of all such functions. For every $\phi\in\mbox{\boldmath$\Phi$}$\ we consider
the adjoint function
\begin{equation}
\phi^{\ast}(t)=t\phi(1/t)\quad\text{where}\quad\phi^{\ast}\in
\mbox{\boldmath$\Phi$},\text{ }(\phi^{\ast})^{\ast}=\phi.\label{0}%
\end{equation}
For every $\phi\in\mbox{\boldmath$\Phi$}$ we consider $\phi$%
\textsl{-divergence} of probability measures $P$ and $Q$ on a measurable space
$({\mathcal{X}},{\mathcal{A}})$ with densities $p,\,q$ w.r.t. a dominating
$\sigma$-finite measure $\lambda.$ In this paper we deal with $P,\,Q$ which
are either measure-theoretically equivalent (i.e. satisfying $pq>0$
\ $\lambda$-a.\thinspace s., in symbols $P\equiv Q$) or measure-theoretically
orthogonal (i.e. satisfying $pq=0$ \ $\lambda$-a.\thinspace s., in symbols
$P\bot Q$).\ Thus, by Liese and Vajda (1987 or 2006), for all $P,Q$ under
consideration
\begin{equation}
D_{\phi}(P,Q)=\left\{
\begin{array}
[c]{ll}%
\int\phi\left(  p/q\right)  \mathrm{d}Q\medskip & \mbox{ \ if \ }\ P\equiv
Q\medskip\\
\phi(0)+\phi^{\ast}(0) & \mbox{ \ if \ }\ P\bot Q
\end{array}
\right. \label{1}%
\end{equation}
where the range of values is%
\begin{equation}
0\leq D_{\phi}(P,Q)\leq\phi(0)+\phi^{\ast}(0)
\end{equation}
and $D_{\phi}(P,Q)=0$\ iff $P=Q$\ or $D_{\phi}(P,Q)=\phi(0)+\phi^{\ast}%
(0)$\ if (for $\phi(0)+\phi^{\ast}(0)<\infty$\ iff) $P\bot Q$. Another
important property is the skew symmetry%
\begin{equation}
D_{\phi}(Q,P)=D_{\phi^{\ast}}(P,Q).\label{1b}%
\end{equation}

We shall deal mainly with the power divergences
\begin{equation}
D_{\alpha}(P,Q):=D_{\phi_{\alpha}}(P,Q)\quad\text{of real powers }\alpha
\in{\mathbb{R}}\medskip\label{5}%
\end{equation}
for the power functions $\phi_{\alpha}\in\mbox{\boldmath$\Phi$}$ defined by
\begin{equation}
\phi_{\alpha}(t)={\frac{t^{\alpha}-\alpha t+\alpha-1}{\alpha(\alpha-1)}}%
\quad\mbox{if}\quad\alpha(\alpha-1)\neq0\label{6}%
\end{equation}
and otherwise by the corresponding limits
\begin{equation}
\phi_{0}(t)=-\ln t+t-1,\qquad\phi_{1}(t)=\phi_{0}^{\ast}(t)=t\ln
t-t+1.\label{7}%
\end{equation}
It is easy to verify for all $\alpha\in{\mathbb{R}}$\ the relation%
\[
\phi_{\alpha}^{\ast}=\phi_{1-\alpha}\text{ so that \ }D_{\alpha}%
(Q,P)=D_{1-\alpha}(P,Q).\medskip
\]
For $P\equiv Q$ we get from (\ref{1}) and (\ref{5})\thinspace--\thinspace
(\ref{7})%
\begin{equation}
D_{\alpha}(P,Q)=\left\{
\begin{array}
[c]{ll}%
\frac{{\Large 1}}{{\Large \alpha(\alpha-1)}}\left[  \int\left(  {p/q}\right)
^{\alpha}\mathrm{d}Q-1\right]  \medskip & \mbox{ \ if \ }\text{ \ \ \ }%
\alpha(\alpha-1)\neq0\medskip\\
\int\ln({p/q)}\,\mathrm{d}P=D_{0}(Q,P) & \mbox{ \ if \ }\text{ \ \ \ }\alpha=1
\end{array}
\right. \label{8}%
\end{equation}
and for $P\bot Q$ similarly%
\begin{equation}
D_{\alpha}(P,Q)=\left\{
\begin{array}
[c]{ll}%
1/\alpha(1-\alpha) & \text{ \ \ \ }\mbox{if}\text{\ \ \ \ }0<\alpha
<1\medskip\\
\infty & \text{ \ \ \ }\mbox{otherwise}.
\end{array}
\right. \label{10}%
\end{equation}
The special cases $D_{2}(P,Q)$ or $D_{1}(P,Q)$ are sometimes called Pearson or
Kullback divergences and $D_{-1}(P,Q)=D_{2}(Q,P)$ or $D_{0}(P,Q)=D_{1}(Q,P)$
reversed Pearson or reverse Kullback divergences, respectively. \medskip

The $\phi$-divergences and power divergences will be applied in the
\textsl{standard statistical estimation model} with i.i.d. observations
$X_{1},\ldots,X_{n}$ governed by $P_{\theta_{0}}$ from a family ${\mathcal{P}%
}=\{P_{\theta}:\theta\in\Theta\}$ of probability measures on $({\mathcal{X}%
},{\mathcal{A}})$ indexed by a set of parameters $\Theta\subset{\mathbb{R}%
}^{d}$. The parameter $\theta_{0}$\ is assumed to be identifiable and the
family ${\mathcal{P}}$ measure-theoreticaly equivalent in the sense%
\begin{equation}
P_{\theta}\neq P_{\theta_{0}}\quad\mbox{and}\quad P_{\theta}\equiv
P_{\theta_{0}}\quad\text{for all}\ \theta,\,\theta_{0}\in\Theta\ \text{\ with
\ }\theta\neq\theta_{0}.\label{12}%
\end{equation}
Further, the family is assumed to be \textsl{continuous} (nonatomic) in the
sense%
\begin{equation}
P_{\theta}(\{x\})=0\text{ \ \ \ for all \ }x\in\mathcal{X},\text{ }\theta
\in\Theta\label{11a}%
\end{equation}
and dominated by a $\sigma$-finite measure $\lambda$\ with densities
\begin{equation}
p_{\theta}=\mathrm{d}P_{\theta}/\mathrm{d}\lambda\quad\text{for all }\theta
\in\Theta.\label{11}%
\end{equation}
In this model the parameter $\theta_{0}$ is assumed to be estimated on the
basis of observations $X_{1},\ldots,X_{n}$ by measurable functions $\theta
_{n}:{\mathcal{X}}^{n}\mapsto\Theta$ called estimates. Collection of estimates
for various sample sizes $n$ is an estimator. Estimators are denoted in this
paper by the same symbols $\theta_{n}$\ as the corresponding
estimates.\medskip

The assumed strict convexity of $\phi(t)$ at $t=1$ together with the
identifiability of $\theta_{0}$\ assumed in (\ref{12}) means that $D_{\phi
}(P_{\theta},P_{\theta_{0}})\geq0$ for all $\theta,\,\theta_{0}\in\Theta
$\ with the equality iff $\theta=\theta_{0}$. In other words, the unknown
parameter $\theta_{0}$ is the unique minimizer of the function $D_{\phi
}(P_{\theta},P_{\theta_{0}})\ $of variable $\theta\in\Theta$,%
\begin{equation}
\theta_{0}=\mathrm{argmin}_{\theta}D(P_{\theta},P_{\theta_{0}})\text{
\ \ \ for every }\theta_{0}\in\Theta.\label{9a}%
\end{equation}
Further, the observations $X_{1},\ldots,X_{n}$ are in a statistically
sufficient manner represented by the empirical probability measure
\begin{equation}
P_{n}={\frac{1}{n}}\sum_{i=1}^{n}P_{X_{i}}\label{14}%
\end{equation}
where $P_{x}$ denotes the Dirac probability measure with all mass concentrated
at $x\in{\mathcal{X}}$. The empirical probability measures $P_{n}$\ are known
to converge weakly to $P_{\theta_{0}}$ as $n\rightarrow\infty$. Therefore by
plugging in (\ref{9a}) the measures $P_{n}$ for $P_{\theta_{0}}$\ one
intuitively expects to obtain the estimator
\begin{equation}
\theta_{n}=\theta_{n,\phi}:=\mbox{{\rm argmin}}_{\theta}D_{\phi}\left(
P_{\theta},P_{n}\right) \label{15}%
\end{equation}
which estimates $\theta_{0}$ consistently in the usual sense of the
convergence $\theta_{n}\rightarrow\theta_{0}$\ for $n\rightarrow\infty$.
However, the reality is different: the problem is that for the continuous
family ${\mathcal{P}}$ under consideration and the discrete family
${\mathcal{P}}_{\text{emp}}$ of empirical distributions (\ref{14}) for which%
\begin{equation}
P_{\theta}\bot P_{n}\Rightarrow D_{\phi}(P_{\theta},P_{n})=\phi(0)+\phi^{\ast
}(0)\quad\mbox{when}\ P_{\theta}\in{\mathcal{P}}\text{ and }P_{n}%
\in{\mathcal{P}}_{\text{emp}}.\label{16}%
\end{equation}
This means that the estimates $\theta_{n}$\ proposed in (\ref{15}) are
trivial, with the $\mathrm{argmin}$ $=\Theta$. \medskip

In the following two sections we list and motivate several modifications of
the minimum divergence rule (\ref{15}) which allow to bypass the problem
(\ref{16}). Some of them are new and some known from the previous literature.
We illustrate the general forms of these estimators by applying them to the
basic standard statistical families and investigate their robustness. The
model of robust statisticians is richer than the standard statistical model
defined by the triplet%
\[
(\mathcal{X},\mathcal{A},\mathcal{Q})\text{ \ \ with \ }\mathcal{Q}%
=\mathcal{P}\cup\mathcal{P}_{\text{emp}}%
\]
introduced above. Namely in addition to the hypothesis that the observations
$X_{1},\ldots,X_{n}$ are i.i.d. by $P_{\theta_{0}}\in{\mathcal{P}}$\ the model
of robust statistics admits the alternative that the observations are
distributed by a probability measure $P_{0}\notin\mathcal{P}$\ with density
\[
\frac{\mathrm{d}P_{0}}{\mathrm{d}\lambda}=p_{0}.
\]

Throughout this paper we assume that $P_{0}$\ is measure-theoretically
equivalent with the probability measures from $\mathcal{P}$\ and we consider
the probability measures%
\begin{equation}
P\in\mathcal{P}\text{ \ and \ }Q\in\mathcal{Q}={\mathcal{P}}^{+}%
\cup{\mathcal{P}}_{\text{emp}}\text{\ \ where \ }{\mathcal{P}}^{+}%
=\mathcal{P}\cup\{P_{0}\}.\label{16b}%
\end{equation}
Measures $P,Q$ are either measure-theoretically equivalent (if $Q\in
{\mathcal{P}}^{+}$) or measure-theoretically orthogonal (if $Q\in{\mathcal{P}%
}_{\text{emp}}$). Therefore the $\phi$-divergences $D_{\phi}(P,Q)$\ are well
defined\ by (\ref{1}) for all pairs $P,Q$\ considered in this paper. Further,
we denote by $\mathbb{L}_{1}(Q)$ the set of all absolutely $Q$-integrable
functions $f:{\mathcal{X\mapsto}}\mathbb{R}$ and put for brevity%
\begin{equation}
Q\cdot f=\int f\,\mathrm{d}Q\text{ \ \ for }f\in\mathbb{L}_{1}(Q).\label{17}%
\end{equation}
\medskip\medskip

In the rest of this section we introduce basic concepts and results of the
robust statistics needed in the sequel. Let us consider\ the Dirac probability
measures $\delta_{x}\in{\mathcal{P}}_{\text{emp}},$ $x\in\mathcal{X}$\ and
denote by \ $C(\mathcal{Q})$\ the set of the convex mixtures%
\begin{equation}
Q_{x,\varepsilon}=(1-\varepsilon)Q+\varepsilon\delta_{x}\text{ \ \ for all
\ }x\in\mathcal{X},\text{\ }Q\in\mathcal{Q}\text{ \ and }0\leq\varepsilon
\leq1.\label{mix}%
\end{equation}
Further, consider a mapping $M(Q,\theta):C(\mathcal{Q})\otimes\Theta
\rightarrow\mathbb{R}\ $differentiable in $\theta\in\Theta$\ for each $Q\in
C(\mathcal{Q})$\ with the derivatives%
\begin{equation}
\Psi(Q,\theta)=\frac{\mathrm{d}}{\mathrm{d}\theta}M(Q,\theta)\label{bigpsi}%
\end{equation}
and let $T(Q)\in\Theta$ solve the equation $\Psi(Q,\theta)=0$ in the variable
$\theta\in\Theta$\ for $Q\in C(\mathcal{Q})$.$\ $The following definition and
theorem deal with the general $M$-estimators
\[
\theta_{n}=\text{\textrm{argmin}}_{\theta}M(P_{n},\theta)\text{ \ \ i.e.
\ \ }\theta_{n}=T(P_{n})\text{ \ for \ }P_{n}\in\mathcal{P}_{\text{emp}%
}\text{.}%
\]
Both the definition and theorem are variants of the well known classical
results of robust statistics, see e.g. Hampel et al. (1986).

\paragraph{Definition 1.1.}

If for some $Q\in{\mathcal{P}}^{+}$ the limits
\begin{equation}
\mbox{IF}(x;T,Q)=\lim_{\varepsilon\downarrow0}{\frac{T(Q_{\varepsilon
,x})-T(Q)}{\varepsilon}}\label{f3}%
\end{equation}
exist for all $x\in{\mathcal{X}}$ then (\ref{f3}) is called influence function
of the estimator $\theta_{n}$ on ${\mathcal{X}}$ at $Q$. \medskip

In the following theorem we consider the functions%
\begin{equation}
\mbox{\boldmath$\psi$}(x,\theta)=\Psi(\delta_{x},\theta)\label{f00}%
\end{equation}
and assume the existence of the derivatives
\begin{equation}
\mathbf{\mathring{\psi}}(x,\theta)=\left(  {\frac{\mathrm{d}}{\mathrm{d}%
\theta}}\right)  ^{\text{t}}\mbox{\boldmath$\psi$}(x,\theta)\quad\text{on
\ }{\mathcal{X}}\otimes\Theta\text{ \ \ ( \ \ with }^{\text{t}}\text{\ for
\textit{transpose})}\label{f5}%
\end{equation}
as well as the expectations%
\begin{equation}
\mbox{\boldmath$I$}(Q)=Q\cdot\mathbf{\mathring{\psi}}(x,T(Q)),\quad
Q\in{\mathcal{P}}^{+}.\label{f6}%
\end{equation}

%%%

\paragraph{Theorem 1.1.}

If the influence function (\ref{f3}) exists then it is given by the formula
\begin{equation}
\mbox{IF}(x;T,Q)=-\mbox{\boldmath$I$}(Q)^{-1}%
\,\mbox{\boldmath$\psi$}(x,T(Q))\label{f7}%
\end{equation}
for the inverse matrix (\ref{f6}).

\paragraph{Proof.}

By definition of $T$, for any $Q\in{\mathcal{P}}^{+}$ and $Q_{\varepsilon,x}%
$considered in (\ref{mix}) it holds
\begin{align*}
0  & ={\frac{Q_{\varepsilon,x}\cdot\mbox{\boldmath$\psi$}(x,T(Q_{\varepsilon
,x}))-Q\cdot\mbox{\boldmath$\psi$}(x,T(Q))}{\varepsilon}}\bigskip\medskip\\
& ={\frac{Q\cdot\lbrack\mbox{\boldmath$\psi$}(x,T(Q_{\varepsilon
,x}))-\mbox{\boldmath$\psi$}(x,T(Q))]}{\varepsilon}}+(\delta_{x}%
-Q)\cdot\mbox{\boldmath$\psi$}(x,T(Q_{\varepsilon,x})).
\end{align*}
Here
\begin{align*}
& \lim_{\varepsilon\downarrow0}{\frac{Q\cdot\lbrack
\mbox{\boldmath$\psi$}(x,T(Q_{\varepsilon,x}%
))-\mbox{\boldmath$\psi$}(x,T(Q))]}{\varepsilon}}\medskip\\
& =Q\cdot\left[  \left(  {\frac{\mathrm{d}}{\mathrm{d}\theta}}\right)
^{\text{t}}\mbox{\boldmath$\psi$}(x,\theta)\right]  _{\theta=T(Q)}%
.\lim_{\varepsilon\downarrow0}{\frac{T(Q_{\varepsilon,x})-T(Q)}{\varepsilon}%
}\medskip\medskip\\
& =Q\cdot\mathbf{\mathring{\psi}}(x,T(Q)).\mbox{IF}(x;T,Q)
\end{align*}
and
\begin{align*}
& \lim_{\varepsilon\downarrow0}(\delta_{x}-Q)\cdot\mathbf{\psi}%
(x,T(Q_{\varepsilon,x}))\bigskip\medskip\\
& =\lim_{\varepsilon\downarrow0}\left[
\mbox{\boldmath$\psi$}(x,T(Q_{\varepsilon,x}))-Q\cdot
\mbox{\boldmath$\psi$}(x,T(Q_{\varepsilon,x}))\right]  \bigskip\medskip\\
& =\mbox{\boldmath$\psi$}(x,T(Q))-Q\cdot
\mbox{\boldmath$\psi$}(x,T(Q))=\mbox{\boldmath$\psi$}(x,T(Q)).
\end{align*}
Therefore we have proved the relation
\[
0=\mbox{\boldmath$I$}(Q).\mbox{IF}(x;T,Q)+\mbox{\boldmath$\psi$}(x,T(Q))
\]
which implies (\ref{f7}).\hfill$\blacksquare$

\bigskip

The estimator $\theta_{n}=T(P_{n})$\ is said to be \textsl{Fisher consistent}
if
\begin{equation}
T(P_{\theta})=\theta\quad\mbox{for all}\ \theta\in\Theta.\label{f4}%
\end{equation}
In the following Corollary and in the sequel, we put%
\begin{equation}
\mbox{IF}(x;T,\theta)=\mbox{IF}(x;T,P_{\theta})\text{ \ \ and \ \ }%
\mbox{\boldmath$I$}(\theta)=\mbox{\boldmath$I$}(P_{\theta})\text{ \ \ (cf.
(\ref{f6})).}\label{f66}%
\end{equation}

\paragraph{Corollary 1.1.}

The influence function of a Fisher consistent estimator at $Q=P_{\theta}$\ is
\begin{equation}
\mbox{IF}(x;T,\theta)=-\mbox{\boldmath$I$}(\theta)^{-1}%
\,\mbox{\boldmath$\psi$}(x,\theta).\label{if}%
\end{equation}

%%%%%%%%%%%%%%%%%%%%%%%%%%%%%%%%%%%%%%%%%%%%%%%%%

\section{SUBDIVERGENCES AND SUPERDIVERGENCES}

Throughout this section we use the likelihood ratios $\boldsymbol{\ell
}_{\theta,\tilde{\theta}}={p_{\theta}}/p_{\tilde{\theta}}$ well defined
a.\thinspace s. on ${\mathcal{X}}$ in the statistical model under
consideration, the nonincreasing functions%
\begin{equation}
\phi^{\#}(t)=\phi(t)-t\phi^{\prime}(t)\quad\text{for every }\phi
\in\mbox{\boldmath$\Phi$}\label{16c}%
\end{equation}
where $\phi^{\prime}$ denotes the derivative of $\phi$, and we restrict
ourselves to the families ${\mathcal{P}}$\ such that%
\begin{equation}
\left\{  \phi\left(  \boldsymbol{\ell}_{\theta,\tilde{\theta}}\right)  ,\text{
}\phi^{\prime}\left(  \boldsymbol{\ell}_{\theta,\tilde{\theta}}\right)
,\text{ }\phi^{\#}\left(  \boldsymbol{\ell}_{\theta,\tilde{\theta}}\right)
\right\}  \subset\mathbb{L}_{1}(Q)\text{ \ \ for all }\theta,\,\tilde{\theta
}\in\Theta\text{ \ and \ }Q\in\mathcal{Q}.\label{16a}%
\end{equation}
Obviously, this assumption automatically holds for all $Q=P_{n}\in
{\mathcal{P}}_{\text{emp}}$. Finally, for all pairs $\theta,\,\tilde{\theta
}\in\Theta$ we consider the functions $L_{\phi}(\theta,\,\tilde{\theta
})=L_{\phi}(\theta,\,\tilde{\theta},\,x)$ of variable $x\in{\mathcal{X}}$
defined by the formula
\[
L_{\phi}(\theta,\,\tilde{\theta})=P_{\theta}\cdot\phi^{\prime}%
(\boldsymbol{\ell}_{\theta,\tilde{\theta}})+\phi^{\#}(\boldsymbol{\ell
}_{\theta,\tilde{\theta}})\text{.}%
\]
$\medskip$

Due to (\ref{16a}), the functions $L_{\phi}(\theta,\,\tilde{\theta})$ are
$Q$-integrable for all $Q\in\mathcal{Q}$. Consider the family of finite
expectations
\begin{equation}
\text{\b{D}}_{\phi,\tilde{\theta}}\left(  P_{\theta},Q\right)  =Q\cdot
L_{\phi}(\,\theta,\tilde{\theta})=P_{\theta}\cdot\phi^{\prime}%
(\boldsymbol{\ell}_{\theta,\tilde{\theta}})+Q\cdot\phi^{\#}(\boldsymbol{\ell
}_{\theta,\tilde{\theta}}),\text{ \ }(P_{\theta},\ Q)\in\mathcal{P}%
\otimes\mathcal{Q}\label{17a}%
\end{equation}
parametrized by $(\phi,\tilde{\theta})\in\mathbf{\Phi\otimes\Theta}$.
Broniatowski \& Keziou (2006) and Liese \& Vajda (2006) independently
established a general supremal representation of $\phi$-divergences $D_{\phi
}\left(  P,Q\right)  $ which implies the following result.

\paragraph{Theorem 2.1.}

For each $(P_{\theta},P_{\theta_{0}})\in\mathcal{P\otimes P}$ and $\phi
\in\mathbf{\Phi}$, the $\phi$-divergence $D_{\phi}\left(  P_{\theta}%
,P_{\theta_{0}}\right)  $ is maximum of the finite expectations \textit{\b{D}%
}$_{\phi,\tilde{\theta}}\left(  P_{\theta},P_{\theta_{0}}\right)  $ over
$\tilde{\theta}\in\Theta$ attained at the unique point $\tilde{\theta}%
=\theta_{0}.$ In other words,%
\begin{equation}
D_{\phi}\left(  P_{\theta},P_{\theta_{0}}\right)  \geq\text{\b{D}}%
_{\phi,\tilde{\theta}}\left(  P_{\theta},P_{\theta_{0}}\right)  \text{ \ \ for
all }\theta,\theta_{0}\in\Theta\label{17f}%
\end{equation}
where the equality holds iff $\tilde{\theta}=\theta_{0}.$

\paragraph{Proof.}

For the sake of completeness we present the simple proof of Liese and Vajda.
For fixed $s>0$, the strictly convex function $\phi(t)$\ is strictly above the
straight line $\phi(s)+\phi^{\prime}(s)(t-s)$ except $t=s$, i.e.%
\[
\phi(t)\geq\phi(s)+\phi^{\prime}(s)(t-s)
\]
with the equality only for $t=s$. Putting in this inequality
$t=\boldsymbol{\ell}_{\theta,\theta_{0}}$, $s=\boldsymbol{\ell}_{\theta
,\tilde{\theta}} $\ and integrating both sides over $P_{\theta_{0}}$\ we get
(\ref{17f}) including the iff condition for the equality.\hfill$\blacksquare
$\bigskip

Theorem 2.1 implies the formula%
\begin{equation}
D_{\phi}\left(  P_{\theta},Q\right)  =\max_{\tilde{\theta}\in\Theta}%
\text{\b{D}}_{\phi,\tilde{\theta}}\left(  P_{\theta},Q\right)  \text{ \ \ for
all }(P_{\theta},Q)\in\mathcal{P}\otimes\mathcal{P}\label{17g}%
\end{equation}
which justifies us to interpret \b{D}$_{\phi,\tilde{\theta}}\left(  P_{\theta
},Q\right)  $\ as $\boldsymbol{subdivergences}$\ of $P_{\theta},Q$ with
parameters $(\phi,\tilde{\theta})\in\mathbf{\Phi\otimes\Theta}.\medskip$

Now we introduce the family of suprema
\begin{equation}
\text{\={D}}_{\phi}\left(  P_{\theta},Q\right)  :=\sup_{\tilde{\theta}%
\in\Theta}\text{\b{D}}_{\phi,\tilde{\theta}}\left(  P_{\theta},Q\right)
\text{ \ \ for all }(P_{\theta},Q)\in\mathcal{P\otimes Q}\label{17e}%
\end{equation}
parametrized by $\phi\in\mathbf{\Phi}$. This family extends the $\phi
$-divergences $D_{\phi}\left(  P,Q\right)  $ from the domain
$\mathcal{P\otimes P}$ to $\mathcal{P\otimes Q}$. Indeed, by Theorem 2.1,
\begin{equation}
\text{\={D}}_{\phi}\left(  P_{\theta},Q\right)  =D_{\phi}\left(  P_{\theta
},Q\right)  \text{ \ for all }(P_{\theta},Q)\in\mathcal{P\otimes P}%
\text{\ .}\label{18}%
\end{equation}
This justifies us to interpret \={D}$_{\phi}\left(  P_{\theta},Q\right)  $ as
$\boldsymbol{superdivergences}$ of $(P_{\theta},Q)\in\mathcal{P\otimes Q}$
with parameters $\phi\in\mathbf{\Phi}$. $\medskip$

Note that (\ref{18})\ need not hold for $Q\notin\mathcal{P}$ because if
$Q=P_{n}\in\mathcal{P}_{\text{emp}}$ then the superdivergence values
\={D}$_{\phi}\left(  P_{\theta},P_{n}\right)  $\ differ from the constant
divergence values $D_{\phi}\left(  P_{\theta},P_{n}\right)  \equiv\phi
(0)+\phi^{\ast}(0)$\ (cf. (\ref{16})). $\medskip$

The subdivergences \b{D}$_{\phi,\tilde{\theta}}\left(  P_{\theta}%
,P_{n}\right)  $\ and superdivergences \={D}$_{\phi}\left(  P_{\theta}%
,P_{n}\right)  $\ can replace the divergences $D_{\phi}\left(  P_{\theta
},P_{n}\right)  $ as optimality criteria in definition of $M$-estimators$.$%
\ Let us consider the families of functionals $\tilde{T}_{\phi,\theta
}:\mathcal{Q}\mapsto\Theta$\ and $T_{\phi}:\mathcal{Q}\ \mapsto\Theta
$\ defined by%
\begin{equation}
\tilde{T}_{\phi,\theta}(Q)=\text{\textrm{argmax}}_{\tilde{\theta}}\text{ \b
{D}}_{\phi,\tilde{\theta}}\left(  P_{\theta},Q\right)  \text{ \ for }%
(\phi,\theta)\in\mathbf{\Phi\otimes\Theta}\label{g1}%
\end{equation}
and%
\begin{equation}
T_{\phi}(Q)=\text{\textrm{argmin}}_{\theta}\text{ \={D}}_{\phi}\left(
P_{\theta},Q\right)  \text{ \ for }\phi\in\mathbf{\Phi}\label{g2}%
\end{equation}
respectively. Replacing the general argument $Q$\ by $P_{n}\,$defined by
(\ref{14})\ we obtain the $\boldsymbol{maximum}$ $\boldsymbol{subdivergence}$
$\boldsymbol{estimators}$ (briefly, the $\max$\b{D}$_{\phi}$-estimators)
\begin{align}
\tilde{\theta}_{\phi,\theta,n}  & \hspace*{-2mm}=\hspace*{-2mm}\tilde{T}%
_{\phi,\theta}(P_{n})=\text{\textrm{argmax}}_{\tilde{\theta}}\text{ \b{D}%
}_{\phi,\tilde{\theta}}\left(  P_{\theta},P_{n}\right)  \medskip\label{G1}\\
& \hspace*{-2mm}=\hspace*{-2mm}\text{\textrm{argmax}}_{\tilde{\theta}}\text{
}\left[  P_{\theta}\cdot\phi^{\prime}(\boldsymbol{\ell}_{\theta,\tilde{\theta
}})+P_{n}\cdot\phi^{\#}(\boldsymbol{\ell}_{\theta,\tilde{\theta}})\right]
\text{ \ \ \ (cf. (\ref{17a}))}\medskip\nonumber\\
& \hspace*{-2mm}=\hspace*{-2mm}\text{\textrm{argmax}}_{\tilde{\theta}}\left[
P_{\theta}\cdot\phi^{\prime}\left(  \frac{p_{\theta}}{p_{\tilde{\theta}}%
}\right)  +\frac{1}{n}\sum_{i=1}^{n}\phi^{\#}\left(  \frac{p_{\theta}(X_{i}%
)}{p_{\tilde{\theta}}(X_{i})}\right)  \right] \label{G1+}%
\end{align}
with escort parameters $\theta\in\mathbf{\Theta,}$ and the
$\boldsymbol{minimum}$ $\boldsymbol{superdivergence}$ $\boldsymbol{estimators}%
$ (briefly, the $\min${\={D}}$_{\phi}${-estimators})
\begin{align}
\theta_{\phi,n}  & \hspace*{-2mm}=\hspace*{-2mm}T_{\phi}(P_{n}%
)=\text{\textrm{argmin}}_{\theta}\text{ \={D}}_{\phi}\left(  P_{\theta}%
,P_{n}\right)  =\text{\textrm{argmin}}_{\theta}\text{\textrm{sup}}%
_{\tilde{\theta}}\text{ \b{D}}_{\phi,\tilde{\theta}}\left(  P_{\theta}%
,P_{n}\right)  \text{\ (cf. (\ref{17e}))}\medskip\label{G2}\\
& \hspace*{-2mm}=\hspace*{-2mm}\text{\textrm{argmin}}_{\theta}%
\text{\textrm{sup}}_{\tilde{\theta}}\left[  P_{\theta}\cdot\phi^{\prime
}(\boldsymbol{\ell}_{\theta,\tilde{\theta}})+P_{n}\,\cdot\phi^{\#}%
(\boldsymbol{\ell}_{\theta,\tilde{\theta}})\right]  \text{
\ \ \ \ \ \ \ \ \ \ \ \ \ \ \ \ \ \ \ \ \ \ (cf. (\ref{17a}))}\medskip
\nonumber\\
& \hspace*{-2mm}=\hspace*{-2mm}\text{\textrm{argmin}}_{\theta}%
\text{\textrm{sup}}_{\tilde{\theta}}\left[  P_{\theta}\cdot\phi^{\prime
}\left(  \frac{p_{\theta}}{p_{\tilde{\theta}}}\right)  +\frac{1}{n}\sum
_{i=1}^{n}\phi^{\#}\left(  \frac{p_{\theta}(X_{i})}{p_{\tilde{\theta}}(X_{i}%
)}\right)  \right]  .\label{G2+}%
\end{align}

\paragraph{Theorem 2.2.}

The $\max$\b{D}$_{\phi}$-estimators are as well as the $\min${\={D}}$_{\phi}%
${-estimators are }Fisher consistent.

\paragraph{Proof.}

By (\ref{17g}) and (\ref{18}),%
\begin{equation}
\tilde{T}_{\phi,\theta}(P_{\theta_{0}})=\text{\textrm{argmax}}_{\tilde{\theta
}}\text{ \b{D}}_{\phi,\tilde{\theta}}\left(  P_{\theta},P_{\theta_{0}}\right)
\text{ \ for }(\phi,\theta)\in\mathbf{\Phi\otimes\Theta}\label{F1}%
\end{equation}
and%
\begin{equation}
T_{\phi}(P_{\theta_{0}})=\text{\textrm{argmin}}_{\theta}\text{ \={D}}_{\phi
}\left(  P_{\theta},P_{\theta_{0}}\right)  \text{ \ for }\phi\in\mathbf{\Phi
}\label{F2}%
\end{equation}
which completes the proof.\hfill$\blacksquare\medskip$

The $\min$\={D}$_{\phi}$-estimators were proposed independently by Liese \&
Vajda (2006) under the name \textbf{modified }$\phi$\textbf{-divergence
estimators}\textsl{\ }and Broniatowski \& Keziou (2006) under the name
\textbf{minimum dual }$\phi$\textbf{-divergence estimators} . The max\b
{D}$_{\phi}$-estimators were proposed by Broniatowski and Keziou (2009) and
called \textbf{dual }$\phi$\textbf{-divergence estimators}\textsl{\ }by
them.\textsl{\ }Both types of these estimators were in the cited papers
motivated by the mentioned Fisher consistency and by the property easily
verifiable from (\ref{G1+}) and (\ref{G2+}), namely that $\phi(t)=-\ln
t$\ implies%
\begin{equation}
\tilde{\theta}_{\phi,\theta,n}=\text{\textrm{argmax}}_{\tilde{\theta}}\text{
}\Sigma_{i=1}^{n}\ln p_{\tilde{\theta}}(X_{i})\text{ \ and \ }\theta_{\phi
,n}=\text{\textrm{argmax}}_{\theta}\text{ }\Sigma_{i=1}^{n}\ln p_{\theta
}(X_{i})\medskip\label{G3}%
\end{equation}
where the left equality holds for all escort parameters $\theta\in\Theta$. In
other words, the logarithmic choice $\phi(t)=-\ln t$ reduces all the variants
of the max\b{D}$_{\phi}$-estimator\ as well as the $\min$\={D}$_{\phi}%
$-estimator to the MLE. It is challenging to investigate the extent to which
the max\b{D}$_{\phi}$-estimators $\tilde{\theta}_{\phi,\theta,n}$\ and the
$\min$\={D}$_{\phi}$-estimator $\theta_{\phi,n}$\ as extensions of the MLE are
efficient and robust under various specifications of $\phi,\theta$\ and $\phi$\ respectively.

$\medskip$

In this paper we restrict ourselves to special subclasses of the power
divergences $D_{\alpha}(P,Q):=D_{\phi_{\alpha}}(P,Q)$\ defined by
(\ref{6})\thinspace--\thinspace(\ref{8}). For the power functions
$\phi_{\alpha}$ from (\ref{6}), (\ref{7}) we get the functions
\begin{equation}
\mathring{\phi}_{\alpha}(t):=t\phi_{\alpha}^{\prime}(t)=\left\{
\begin{array}
[c]{ll}%
\frac{{\Large t}^{{\Large \alpha}}{\Large -t}}{{\Large \alpha-1}} &
\mbox{\ \ for }\alpha\neq1\medskip\\
\lim_{\alpha\rightarrow1}\frac{{\Large t}^{{\Large \alpha}}{\Large -t}%
}{{\Large \alpha-1}}=t\ln t\medskip & \mbox{\ \ for }\alpha=1
\end{array}
\right. \label{f1}%
\end{equation}
and
\begin{equation}
\phi_{\alpha}^{\#}(t)=\phi_{\alpha}(t)-\mathring{\phi}_{\alpha}(t)=\left\{
\begin{array}
[c]{ll}%
\frac{{\Large 1}}{{\Large \alpha}}\left(  1-t^{\alpha}\right)  & \mbox{\ \
for }\alpha\neq0\medskip\\
\lim_{\alpha\rightarrow0}\frac{{\Large 1}}{{\Large \alpha}}\left(
1-t^{\alpha}\right)  =-\ln t\medskip & \mbox{\ \ for }\alpha=0.
\end{array}
\right. \label{f2}%
\end{equation}
They lead to the max\b{D}$_{\alpha}$-estimators (briefly, $\boldsymbol{power}$
$\boldsymbol{subdivergence}$ $\boldsymbol{estimators}$)%
\begin{equation}
\tilde{\theta}_{\alpha,\theta,n}=\text{\textrm{argmax}}_{\tilde{\theta}%
}\left[  P_{\tilde{\theta}}\cdot\mathring{\phi}_{\alpha}\left(  \frac
{p_{\theta}}{p_{\tilde{\theta}}}\right)  +P_{n}\cdot\phi_{\alpha}^{\#}\left(
\frac{p_{\theta}}{p_{\tilde{\theta}}}\right)  \right] \label{ggg1}%
\end{equation}
\ with power parameters $\alpha\in\mathbb{R}$\ and escort parameters
$\theta\in\Theta$ and to the $\min$\={D}$_{\alpha}$-estimators\ (briefly,
$\boldsymbol{power}$ $\boldsymbol{superdivergence}$ $\boldsymbol{estimators}$)%
\begin{equation}
\theta_{\alpha,n}=\text{\textrm{argmin}}_{\theta}\text{\textrm{sup}}%
_{\tilde{\theta}}\left[  P_{\tilde{\theta}}\cdot\mathring{\phi}_{\alpha
}\left(  \frac{p_{\theta}}{p_{\tilde{\theta}}}\right)  +P_{n}\cdot\phi
_{\alpha}^{\#}\left(  \frac{p_{\theta}}{p_{\tilde{\theta}}}\right)  \right]
\label{ggg2}%
\end{equation}
with power parameters $\alpha\in\mathbb{R}.$ If the argmaxima in (\ref{ggg1})
exist then%
\begin{equation}
\theta_{\alpha,n}=\text{\textrm{argmin}}_{\theta}\left[  P_{\tilde{\theta
}_{\alpha,\theta,n}}\cdot\mathring{\phi}_{\alpha}\left(  \frac{p_{\theta}%
}{p_{\tilde{\theta}_{\alpha,\theta,n}}}\right)  +P_{n}\cdot\phi_{\alpha}%
^{\#}\left(  \frac{p_{\theta}}{p_{\tilde{\theta}_{\alpha,\theta,n}}}\right)
\right]  .\label{ggg3}%
\end{equation}
$\medskip$

The next two subsections deal correspondingly with the max\b{D}$_{\alpha}%
$-estimators\ and $\min$\={D}$_{\alpha}$-estimators. In both sections are
considered the power parameters $\alpha\geq0$. Since $\phi_{0}(t)=-\ln t,$\ we
see from (\ref{G3}) that%
\begin{equation}
\tilde{\theta}_{0,\theta,n}=\text{\textrm{argmax}}_{\tilde{\theta}}\text{
}\Sigma_{i=1}^{n}\ln p_{\tilde{\theta}}(X_{i})\text{ \ and \ }\theta
_{0,n}=\text{\textrm{argmax}}_{\theta}\text{ }\Sigma_{i=1}^{n}\ln p_{\theta
}(X_{i})\medskip\label{MLE}%
\end{equation}
are the MLE's. If $\alpha>0$ then by (\ref{f1}) - (\ref{ggg2}),%
\begin{equation}
\tilde{\theta}_{\alpha,\theta,n}=\text{\textrm{argmin}}_{\tilde{\theta}}\text{
}M_{\alpha,\theta}(P_{n},\tilde{\theta})\label{gg1}%
\end{equation}
and%
\begin{equation}
\theta_{\alpha,n}=\text{\textrm{argmax}}_{\theta}\text{\textrm{inf}}%
_{\tilde{\theta}}\text{ }M_{\alpha,\theta}(P_{n},\tilde{\theta})\equiv
\text{\textrm{argmax}}_{\theta}\text{ }M_{\alpha,\theta}(P_{n},\tilde{\theta
}_{\alpha,\theta,n})\label{gg2}%
\end{equation}
where%
\begin{align}
M_{\alpha,\theta}(Q,\tilde{\theta})  & =\frac{1}{1-\alpha}P_{\tilde{\theta}%
}\cdot\left(  \frac{p_{\theta}}{p_{\tilde{\theta}}}\right)  ^{\alpha}+\frac
{1}{\alpha}Q\cdot\left(  \frac{p_{\theta}}{p_{\tilde{\theta}}}\right)
^{\alpha}\text{ \ \ if }\alpha>0\text{, }\alpha\neq1\nonumber\\
& \label{gg3}\\
& =P_{\theta}\cdot\ln\frac{p_{\tilde{\theta}}}{p_{\theta}}+Q\cdot
\frac{p_{\theta}}{p_{\tilde{\theta}}}\text{
\ \ \ \ \ \ \ \ \ \ \ \ \ \ \ \ \ \ \ \ \ \ if }\alpha=1\nonumber
\end{align}
for all $Q\in\mathcal{Q}$. $\medskip$

Throughout both subsections we restrict ourselves to the densities $p_{\theta
}$\ twice differentiable with respect to $\theta\in\Theta\subset\mathbb{R}%
^{d}$, we put
\begin{equation}
s_{\theta}=\frac{\mathrm{d}}{\mathrm{d}\theta}\ln p_{\theta}\text{ \ \ and
\ \ }\mathring{s}_{\theta}=\left(  \frac{\mathrm{d}}{\mathrm{d}\theta}\right)
^{\text{t}}s_{\theta}\label{p1}%
\end{equation}
and suppose that the functions $M_{\alpha,\theta}(Q,\tilde{\theta})$\ of
(\ref{gg3}) are twice differentiable in the vector variable $\tilde{\theta}%
,$\ with the differentiation and integration interchangeable in (\ref{gg3}).
Moreover, we suppose that the derivatives%
\begin{equation}
\Psi_{\alpha,\theta}(Q,\tilde{\theta})=\frac{\mathrm{d}}{\mathrm{d}%
\tilde{\theta}}M_{\alpha,\theta}(Q,\tilde{\theta})=P_{\tilde{\theta}}%
\cdot\left(  \frac{p_{\theta}}{p_{\tilde{\theta}}}\right)  ^{\alpha}%
s_{\tilde{\theta}}-Q\cdot\left(  \frac{p_{\theta}}{p_{\tilde{\theta}}}\right)
^{\alpha}s_{\tilde{\theta}}.\label{p3}%
\end{equation}
admit solutions of the equations $\Psi_{\alpha,\theta}(Q,\tilde{\theta})=0
$\ in the variable $\tilde{\theta}\in\Theta$\ for $Q\in\emph{Q}$.

\subsection{Power subdivergence estimators}

In this subsection we study the $\max$\b{D}$_{\alpha}$-estimators
$\tilde{\theta}_{\alpha,\theta,n}$ with the divergence power parameters
$\alpha\geq0$\ and the escort parameters $\theta\in\Theta$. As said above, for
$\alpha=0$\ they coincide with the MLE's (\ref{MLE}). Therefore we restrict
ourselves to\ $\alpha>0$ and to the definition formula (\ref{gg1}),
(\ref{gg3}). $\medskip$

By assumptions, the argminima
\begin{equation}
\tilde{T}_{\alpha,\theta}(Q)=\text{\textrm{argmin}}_{\tilde{\theta}}\text{
}M_{\alpha,\theta}(Q,\tilde{\theta}),\text{ \ \ }\alpha>0,\text{\ \ }%
Q\in\mathcal{Q}\text{\ \ \ \ \ \ \ (cf. (\ref{g1}))}\label{p2}%
\end{equation}
solve the equations $\Psi_{\alpha,\theta}(Q,\tilde{\theta})=0$\ in the
variable $\tilde{\theta}\in\Theta$\ and, in particular, $\tilde{\theta
}_{\alpha,\theta,n}=\tilde{T}_{\alpha,\theta}(P_{n})$\ are for all $\alpha
>0$\ solutions of the equations%
\begin{equation}
P_{\tilde{\theta}}\cdot\left(  \frac{p_{\theta}}{p_{\tilde{\theta}}}\right)
^{\alpha}s_{\tilde{\theta}}-\frac{1}{n}\sum_{i=1}^{n}\left(  \frac{p_{\theta
}(X_{i})}{p_{\tilde{\theta}}(X_{i})}\right)  ^{\alpha}s_{\tilde{\theta}}%
(X_{i})=0\label{p3b}%
\end{equation}
in the variable $\tilde{\theta}\in\Theta$.$\medskip$

\paragraph{Theorem 2.1.1.}

The influence functions of the $\max$\b{D}$_{\alpha}$-estimators
$\tilde{\theta}_{\alpha,\theta,n}$ under consideration are at $P_{\theta_{0}}%
$\ given by the formula%
\begin{align}
\text{IF}(x;\tilde{T}_{\alpha,\theta},\theta_{0})  & =\boldsymbol{I}%
_{\alpha,\theta}(\theta_{0})^{-1}\left[  \left(  \frac{p_{\theta}%
(x)}{p_{\theta_{0}}(x)}\right)  ^{\alpha}s_{\theta_{0}}(x)-P_{\theta_{0}}%
\cdot\left(  \frac{p_{\theta}}{p_{\theta_{0}}}\right)  ^{\alpha}s_{\theta_{0}%
}\right]  \text{ \ if \ }\alpha>0\label{t1}\\
& \nonumber\\
\text{IF}(x;\tilde{T}_{0,\theta},\theta_{0})  & =\boldsymbol{I}(\theta
_{0})^{-1}s_{\theta_{0}}(x)\text{
\ \ \ \ \ \ \ \ \ \ \ \ \ \ \ \ \ \ \ \ \ \ \ \ \ \ \ \ \ \ \ \ \ \ \ \ \ \ \ \ \ \ \ \ \ \ \ otherwise}%
\label{t2}%
\end{align}
where%
\begin{align}
\boldsymbol{I}_{\alpha,\theta}(\theta_{0})  & =P_{\theta_{0}}\cdot\left(
\frac{p_{\theta}}{p_{\theta_{0}}}\right)  ^{\alpha}s_{\theta_{0}}^{\text{t}%
}s_{\theta_{0}}\text{ \ \ \ \ \ if \ }\alpha>0\medskip\medskip\medskip
\label{t3}\\
\boldsymbol{I}(\theta_{0})  & =P_{\theta_{0}}\cdot s_{\theta_{0}}^{\text{t}%
}s_{\theta_{0}}\text{ \ \ \ \ \ \ \ \ \ \ \ \ \ \ \ \ if \ }\alpha
=0\text{.}\label{t4}%
\end{align}
If the escort parameter $\theta$\ coincides with the true parameter
$\theta_{0}$ then%
\[
\text{IF}(x;\tilde{T}_{\alpha,\theta_{0}},\theta_{0})=\boldsymbol{I}%
(\theta_{0})^{-1}s_{\theta_{0}}(x)\text{ \ \ for all }\alpha\geq0\text{.}%
\]

\paragraph{Proof.}

By (\ref{f00}) and (\ref{p3}),%
\begin{equation}
\mathbf{\psi}_{\alpha,\theta}(x,\tilde{\theta})=\Psi_{\alpha,\theta}%
(\delta_{x},\tilde{\theta})=P_{\tilde{\theta}}\cdot\left(  \frac{p_{\theta}%
}{p_{\tilde{\theta}}}\right)  ^{\alpha}s_{\tilde{\theta}}-\delta_{x}%
\cdot\left(  \frac{p_{\theta}}{p_{\tilde{\theta}}}\right)  ^{\alpha}%
s_{\tilde{\theta}}\label{p3a}%
\end{equation}
and under the assumptions stated above%
\begin{equation}
\mathbf{\mathring{\psi}}_{\alpha,\theta}(x,\tilde{\theta})=\left(
\frac{\mathrm{d}}{\mathrm{d}\tilde{\theta}}\right)  ^{\text{t}}\mathbf{\psi
}_{\alpha,\theta}(x,\tilde{\theta})=P_{\tilde{\theta}}\cdot\left(
\frac{p_{\theta}}{p_{\tilde{\theta}}}\right)  ^{\alpha}s^{\text{t}}%
s_{\tilde{\theta}}-P_{\tilde{\theta}}\cdot\Lambda_{\alpha,\theta,\tilde
{\theta}}+\Lambda_{\alpha,\theta,\tilde{\theta}}(x)\label{p33}%
\end{equation}
for%
\[
\Lambda_{\alpha,\theta,\tilde{\theta}}(x)=\left(  \frac{p_{\theta}%
(x)}{p_{\tilde{\theta}}(x)}\right)  ^{\alpha}\left[  \alpha s_{\tilde{\theta}%
}(x)^{\text{t}}s_{\tilde{\theta}}(x)-\mathring{s}_{\tilde{\theta}}(x)\right]
.
\]
Further, by (\ref{f66}), (\ref{f6}) and (\ref{p33}),%
\[
\boldsymbol{I}_{\alpha,\theta}(\theta_{0})=P_{\theta_{0}}\cdot
\mathbf{\mathring{\psi}}_{\alpha,\theta}(x,\theta_{0})=P_{\theta_{0}}%
\cdot\left(  \frac{p_{\theta}}{p_{\theta_{0}}}\right)  ^{\alpha}s^{\text{t}%
}s_{\theta_{0}}\medskip
\]
and (\ref{if}) leads to the\textsl{\ }influence functions%
\[
\text{IF}(x;\tilde{T}_{\alpha,\theta},\theta_{0})=-\boldsymbol{I}%
_{\alpha,\theta}(\theta_{0})^{-1}\mathbf{\psi}_{\alpha,\theta}(x,\theta
_{0}).\medskip
\]
The substitution from (\ref{p3a}) yields the desired formula (\ref{t1}). In
the MLE case $\alpha=0$\ we get for all escort parameters $\theta$\ the
classical MLE influence function (\ref{t2}) with the classical Fisher
information matrix given in (\ref{t4}). This influence function is obtained
also\ if the escort parameter $\theta$\ coincides with the true parameter
$\theta_{0}$ as in this case the estimators with all power parameters
$\alpha\geq0$\ reduce to the MLE (cf. (\ref{MLE})).\hfill$\blacksquare
\medskip$

Next follow special examples of the influence functions (\ref{t1}),
(\ref{t2}).$\medskip$

\paragraph{Example 2.1.1: Power subdivergence estimators in normal family.}

Let the observation space $({\mathcal{X}},{\mathcal{A}})$\ be the Borel line
$(\mathbb{R},{\mathcal{B}})$\ and ${\mathcal{P}}=\{P_{\mu,\sigma}:\mu
\in\mathbb{R},$ $\sigma>0\}$ the normal family with parameters of location
$\mu$\ and scale $\sigma$\ (i.e. variances $\sigma^{2}$). We are interested in
the max\b{D}$_{\alpha}$-estimates $(\tilde{\mu}_{\alpha,\mu,\sigma,n}%
,\tilde{\sigma}_{\alpha,\mu,\sigma,n})$ with power parameters $\alpha\geq0$
and escort parameters $(\mu,\sigma)\in\mathbb{R}\otimes(0,\infty)\}.$
$\medskip$

If $\alpha=0$\ then these estimators\ reduce\ for all escort parameters
$\mu,\sigma$ to the well known MLE's%
\begin{equation}
(\tilde{\mu}_{0,\mu,\sigma,n},\tilde{\sigma}_{0,\mu,\sigma,n})=\left(
\frac{1}{n}\sum_{i=1}^{n}X_{i},\ \sqrt{\frac{1}{n}\sum_{i=1}^{n}\left(
X_{i}-\tilde{\mu}_{0,n}\right)  ^{2}}\right) \label{mle}%
\end{equation}
For $0<\alpha<1$\ the function (\ref{gg3}) takes on the form%
\begin{equation}
M_{\alpha,\mu,\sigma}(Q,\tilde{\mu},\tilde{\sigma})=\frac{1}{1-\alpha
}P_{\tilde{\mu},\tilde{\sigma}}\cdot\left(  \frac{p_{\mu,\sigma}}%
{p_{\tilde{\mu},\tilde{\sigma}}}\right)  ^{\alpha}+\frac{1}{\alpha}%
Q\cdot\left(  \frac{p_{\mu,\sigma}}{p_{\tilde{\mu},\tilde{\sigma}}}\right)
^{\alpha}\label{M}%
\end{equation}
\ where%
\begin{equation}
\left(  \frac{p_{\mu,\sigma}(x)}{p_{\tilde{\mu},\tilde{\sigma}}(x)}\right)
^{\alpha}=\left(  \frac{\tilde{\sigma}}{\sigma}\right)  ^{\alpha}\exp\left\{
\frac{\alpha\left(  x-\tilde{\mu}\right)  ^{2}}{2\tilde{\sigma}^{2}}%
-\frac{\alpha\left(  x-\mu\right)  ^{2}}{2\sigma^{2}}\right\}  ,\label{p7}%
\end{equation}
and%
\begin{equation}
P_{\tilde{\mu},\tilde{\sigma}}\cdot\left(  \frac{p_{\mu,\sigma}}{p_{\tilde
{\mu},\tilde{\sigma}}}\right)  ^{\alpha}=\exp\left\{  -\frac{\alpha
(1-\alpha)(\mu-\tilde{\mu})^{2}}{2[\alpha\tilde{\sigma}^{2}+(1-\alpha
)\sigma^{2}]}-\ln\frac{\sqrt{\alpha\tilde{\sigma}^{2}+(1-\alpha)\sigma^{2}}%
}{\tilde{\sigma}^{\alpha}\sigma^{1-\alpha}}\right\}  .\label{p6}%
\end{equation}
Using the likelihood ratio function (\ref{p7}) and the score function%
\begin{equation}
s_{\mu,\sigma}(x)=\left(  \frac{x-\mu}{\sigma^{2}},\frac{1}{\sigma}\left[
\left(  \frac{x-\mu}{\sigma}\right)  ^{2}-1\right]  \right) \label{p8}%
\end{equation}
one obtains for all $\alpha>0$\ the derivative%
\begin{equation}
\Psi_{\alpha,\mu,\sigma}(Q,\tilde{\mu},\tilde{\sigma})=\left(  \frac
{\mathrm{d}}{\mathrm{d}\tilde{\mu}},\frac{\mathrm{d}}{\mathrm{d}\tilde{\sigma
}}\right)  \text{ }M_{\alpha,\mu,\sigma}(Q,\tilde{\mu},\tilde{\sigma
})=P_{\tilde{\mu},\tilde{\sigma}}\cdot\left(  \frac{p_{\mu,\sigma}}%
{p_{\tilde{\mu},\tilde{\sigma}}}\right)  ^{\alpha}s_{\tilde{\mu},\tilde
{\sigma}}-Q\cdot\left(  \frac{p_{\mu,\sigma}}{p_{\tilde{\mu},\tilde{\sigma}}%
}\right)  ^{\alpha}s_{\tilde{\mu},\tilde{\sigma}}\label{p13}%
\end{equation}
and the max\b{D}$_{\alpha}$-estimators as the argminima%
\begin{equation}
(\tilde{\mu}_{\alpha,\mu,\sigma,n},\tilde{\sigma}_{\alpha,\mu,\sigma
,n})=\text{\textrm{argmin}}_{\tilde{\mu},\tilde{\sigma}}\left[  \frac
{1}{1-\alpha}P_{\tilde{\mu},\tilde{\sigma}}\cdot\left(  \frac{p_{\mu,\sigma}%
}{p_{\tilde{\mu},\tilde{\sigma}}}\right)  ^{\alpha}+\frac{1}{\alpha n}%
\sum_{i=1}^{n}\left(  \frac{p_{\mu,\sigma}(X_{i})}{p_{\tilde{\mu}%
,\tilde{\sigma}}(X_{i})}\right)  ^{\alpha}\right] \label{p15}%
\end{equation}
or, equivalently, as solutions of the equations%
\begin{equation}
P_{\tilde{\mu},\tilde{\sigma}}\cdot\left(  \frac{p_{\mu,\sigma}}{p_{\tilde
{\mu},\tilde{\sigma}}}\right)  ^{\alpha}s_{\tilde{\mu},\tilde{\sigma}}%
-\frac{1}{n}\sum_{i=1}^{n}\left(  \frac{p_{\mu,\sigma}(X_{i})}{p_{\tilde{\mu
},\tilde{\sigma}}(X_{i})}\right)  ^{\alpha}s_{\tilde{\mu},\tilde{\sigma}%
}(X_{i})=0.\label{p16}%
\end{equation}

By Theorem 2.1.1, the influence functions of these estimators at $P_{\mu
_{0},\sigma_{0}}$\ are%
\begin{equation}
\text{IF}(x;\tilde{T}_{\alpha,\mu,\sigma},\mu_{0},\sigma_{0})=\boldsymbol{I}%
_{\mu,\sigma}(\mu_{0},\sigma_{0})^{-1}\left[  \left(  \frac{p_{\mu,\sigma}%
(x)}{p_{\mu_{0},\sigma_{0}}(x)}\right)  ^{\alpha}s_{\mu_{0},\sigma_{0}%
}(x)-P_{\mu_{0},\sigma_{0}}\cdot\left(  \frac{p_{\mu,\sigma}}{p_{\mu
_{0},\sigma_{0}}}\right)  ^{\alpha}s_{\mu_{0},\sigma_{0}}\right] \label{p17}%
\end{equation}
for%
\begin{equation}
\boldsymbol{I}_{\mu,\sigma}(\mu_{0},\sigma_{0})=P_{\mu_{0},\sigma_{0}}%
\cdot\left(  \frac{p_{\mu,\sigma}}{p_{\mu_{0},\sigma_{0}}}\right)  ^{\alpha
}s_{\mu_{0},\sigma_{0}}^{\text{t}}s_{\mu_{0},\sigma_{0}}.\label{p18}%
\end{equation}

\paragraph{Example 2.1.2: Power subdivergence estimators of location.}

Let in the frame of previous example ${\mathcal{P}}=\{P_{\mu}:\mu\in
\mathbb{R}\}$ be the standard normal family with the location parameter $\mu
$\ and scale $\sigma=1$. Then the function (\ref{M}) takes on the form%
\begin{equation}
M_{\alpha,\mu}(Q,\tilde{\mu})=\frac{1}{1-\alpha}\left(  \eta_{\alpha,\mu}%
(\mu,\tilde{\mu})\right)  ^{\alpha-1}+\frac{1}{\alpha}Q\cdot\eta_{\alpha,\mu
}(x,\tilde{\mu})\label{p10}%
\end{equation}
for $\alpha>0,\alpha\neq1$\ where%
\[
\eta_{\alpha,\mu}(x,\tilde{\mu})=\exp\left\{  \alpha(\tilde{\mu}-\mu
)(\tilde{\mu}+\mu-2x)/2\right\}  ,\text{ \ \ }x\in\mathbb{R}.
\]
The max\b{D}$_{\alpha}$-estimates $\tilde{\mu}_{\alpha,\mu,n}$ of location
$\mu_{0}$ with the divergence parameters $0\leq\alpha<1$ and escort parameters
$\mu\in\mathbb{R}$ are the MLE's%
\begin{equation}
\tilde{\mu}_{0,\mu,n}=\boldsymbol{\bar{X}}_{n}=\frac{1}{n}\sum_{i=1}^{n}%
X_{i}\label{63a}%
\end{equation}
if $\alpha=0.$ Otherwise they are the minimizers%
\begin{equation}
\tilde{\mu}_{\alpha,\mu,n}=\text{\textrm{argmin}}_{\tilde{\mu}}M_{\alpha,\mu
}(P_{n},\tilde{\mu})\medskip\label{p9}%
\end{equation}
or, equivalently, solutions of the equations%
\[
\Psi_{\alpha,\mu}(P_{n},\tilde{\mu})=0
\]
in they variable $\tilde{\mu}\in\mathbb{R}$\ for%
\begin{align}
\Psi_{\alpha,\mu}(Q,\tilde{\mu})  & =\frac{\mathrm{d}}{\mathrm{d}\tilde{\mu}%
}M_{\alpha,\mu}(Q,\tilde{\mu})\nonumber\\
& \nonumber\\
& =Q\cdot(\tilde{\mu}-x)\eta_{\alpha,\mu}(x,\tilde{\mu})-\alpha(\tilde{\mu
}-\mu)\eta_{\alpha,\mu}^{\alpha-1}(\mu,\tilde{\mu}).\label{locc}%
\end{align}
$\medskip$

\ Let $\tilde{T}_{\alpha,\mu}(Q)$ be the solution of the equation
$\Psi_{\alpha,\mu}(Q,\tilde{\mu})=0$ in the variable $\tilde{\mu}\in
\mathbb{R}$ and let $Q_{\mu_{0}}$\ denote the shift of the distribution
$Q$\ by $\mu_{0}.$\ Then%
\[
Q_{\mu_{0}}\cdot(\tilde{\mu}-x)\eta_{\alpha,\mu}(x,\tilde{\mu})=Q\cdot
(\tilde{\mu}-\mu_{0}-x)\eta_{\alpha,\mu-\mu_{0}}(x,\tilde{\mu}-\mu
_{0}))\medskip
\]
so that $\tilde{T}_{\alpha,\mu}(Q_{\mu_{0}})=\mu_{0}+\tilde{T}_{\alpha,\mu
-\mu_{0}}(Q).$\ This means that the estimators (\ref{p9}) are Fisher
consistent in the normal family $\mathcal{P}_{\sigma}=\left\{  P_{\mu
_{0},\sigma}=N(\mu_{0},\sigma^{2}):\mu_{0}\in\mathbb{R}\right\}  $\ with
$\sigma>0$ fixed if and only if the solution $\tilde{T}_{\alpha,\mu
}(P_{0,\sigma})$ of the equation%
\begin{equation}
P_{0,\sigma}\cdot(\tilde{\mu}-x)\eta_{\alpha,\mu}(x,\tilde{\mu})-\alpha
(\tilde{\mu}-\mu)\eta_{\alpha,\mu}^{\alpha-1}(\mu,\tilde{\mu})=0\label{49}%
\end{equation}
in the variable $\tilde{\mu}$ satisfies the condition\
\begin{equation}
\tilde{T}_{\alpha,\mu}(P_{0,\sigma})=0\text{ \ \ for all }\mu\in
\mathbb{R}.\medskip\label{p11}%
\end{equation}
By evaluating the function $P_{0,\sigma}\cdot(\tilde{\mu}-x)\eta_{\alpha,\mu
}(x,\tilde{\mu})$ of variables $\sigma,\mu,\tilde{\mu}$ and inserting it in
(\ref{49}), one can verify that (\ref{p11}) holds if and only if $\sigma=1$.
The \textquotedblleft if" part follows from the Fisher consistency of
$\tilde{T}_{\alpha,\mu}$ established in Theorem 2.2 which implies%
\[
\tilde{T}_{\alpha,\mu}(P_{0,1})\equiv\tilde{T}_{\alpha,\mu}(P_{0})=0\text{
\ \ for }P_{0,1}\equiv P_{0}\in\mathcal{P}\text{ and all }\mu\in\mathbb{R.}%
\]
However, the \textquotedblleft only if" assertion is \textbf{new and
surprising} in the sense that it indicates a relatively easy loss of
consistency of the max\b{D}$_{\alpha}$-estimators.\ 

\paragraph{Problem 2.1.1.}

It remains to be verified analytically or by simulations whether the
estimators $\tilde{\mu}_{\alpha,\bar{X}_{n},n}$ with the adaptive MLE escort
parameters $\bar{X}_{n}$ are Fisher consistent under all hypothetical models
$P_{\mu,\sigma}=N(\mu,\sigma^{2}),$\ $\sigma>0$ or, more generally, whether
the adaptive estimators%
\begin{equation}
\tilde{\theta}_{\alpha,\tau_{n},n}\text{ \ \ with the MLE escorts }\tau
_{n}=\tilde{\theta}_{0,n}\text{\ given by (\ref{G3})}\label{G4}%
\end{equation}
are Fisher consistent under the hypothetical models $P_{\theta_{0}}$, and
eventually consistent and robust under contaminated versions of these
models.$\medskip\medskip$

Let us turn to the influence curves IF$(x;T_{\alpha,\mu},\mathbb{\mu}%
_{0}),0<\alpha<1$ at the data source $P_{\mathbb{\mu}_{0}}$. Here
$s_{\mathbb{\mu}_{0}}^{\text{t}}(x)s_{\mathbb{\mu}_{0}}(x)=s_{\mathbb{\mu}%
_{0}}^{2}(x)=(\mathbb{\mu}_{0}-x)^{2}$ so that, by (\ref{f66}) and
(\ref{p18}),%
\begin{align}
I_{\alpha,\mu}(\mathbb{\mu}_{0})  & =\boldsymbol{I}_{\alpha,\mu}%
(P_{\mathbb{\mu}_{0}})=P_{\mathbb{\mu}_{0}}\cdot\left(  \frac{p_{\mu}%
}{p_{\mathbb{\mu}_{0}}}\right)  ^{\alpha}s_{\mathbb{\mu}_{0}}^{2}%
\medskip\nonumber\\
& =\frac{1}{\sqrt{2\pi}}\int(\mathbb{\mu}_{0}-x)^{2}\exp\left\{  -\frac
{\alpha(x-\mu)^{2}+(1-\alpha)(x-\mathbb{\mu}_{0})^{2}}{2}\right\}
\mathrm{d}x\medskip\label{loc3}\\
& =\left[  1+\alpha^{2}(\mathbb{\mu}_{0}-\mu)^{2}\right]  \exp\left\{
\frac{\alpha(\alpha-1)(\mathbb{\mu}_{0}-\mu)^{2}}{2}\right\}  .\nonumber
\end{align}
If we put%
\[
\mathbf{\psi}_{\alpha,\mu}(x,\mathbb{\mu}_{0})=\Psi_{\alpha,\mu}(\delta
_{x},\mathbb{\mu}_{0})=(\mathbb{\mu}_{0}-x)\eta_{\alpha,\mu}(x,\mathbb{\mu
}_{0})-\alpha(\mathbb{\mu}_{0}-\mu)\eta_{\alpha,\mu}^{\alpha-1}(\mu
,\mathbb{\mu}_{0})\text{ \ \ \ (cf. (\ref{locc}))}%
\]
then, by (\ref{p17}),
\begin{align}
\text{IF}(x;T_{\alpha,\mu},\mathbb{\mu}_{0})  & =-\frac{\mathbf{\psi}%
_{\alpha,\mu}(x,\mathbb{\mu}_{0})}{I_{\alpha,\mu}(\mathbb{\mu}_{0}%
)}\nonumber\\
& =\frac{(x-\mathbb{\mu}_{0})e^{\alpha(\mathbb{\mu}_{0}-\mu)(\mathbb{\mu}%
_{0}+\mu-2x)/2}+\alpha(\mathbb{\mu}_{0}-\mu)e^{\alpha(\alpha-1)(\mathbb{\mu
}_{0}-\mu)^{2}/2}}{\left[  1+\alpha^{2}(\mathbb{\mu}_{0}-\mu)^{2}\right]
e^{\alpha(\alpha-1)(\mathbb{\mu}_{0}-\mu)^{2}/2}}.\medskip\label{loc4}%
\end{align}
This formula remains valid also for $\alpha=0$\ because then it reduces to the
well known influence function%
\[
\text{IF}(x;MLE,\mathbb{\mu}_{0})=x-\mathbb{\mu}_{0}%
\]
of the $MLE=T_{0,\mu}$\ which is not depending on the escort parameter $\mu$.
We see that the influence curve (\ref{loc4}) is unbounded for all $\mu,\mu
_{0}\in\mathbb{R}$ and $0\leq\alpha<1$. For $0<\alpha<1$ and the escort
parameters $\mu$\ different from the true $\mu_{0}$\ the influence functions
IF$(x;T_{\alpha,\mu},\mathbb{\mu}_{0})$ contain the constant terms
IF$(\mathbb{\mu}_{0};T_{\alpha,\mu},\mathbb{\mu}_{0})\neq0$\ and, moreover,
increase to infinity exponentially for $x\rightarrow\infty$ or $x\rightarrow
-\infty.$ Therefore $T_{\alpha,\mu}$\ are strongly non-robust.

\paragraph{Example 2.1.3: Power subdivergence estimators of scale.}

Let in the frame of Example 2.1.1, ${\mathcal{P}}=\{P_{\sigma}:\sigma>0\}$ be
the standard normal family with the location parameter $\mu=0$ and scale
$\sigma$ and let us consider the max\b{D}$_{\alpha}$-estimators $\tilde
{\sigma}_{\alpha,\sigma,n}$ of scale $\sigma_{0}$ with the divergence
parameters $0\leq\alpha<1$ and escort parameters $\sigma>0.$ For $\alpha=0$
they reduce to the standard deviations%
\[
\tilde{\sigma}_{0,\sigma,n}=\left(  \frac{1}{n}\sum_{i=1}^{n}\left(
X_{i}-\boldsymbol{\bar{X}}_{n}\right)  ^{2}\right)  ^{1/2}%
\]
and otherwise they are of the form
\[
\tilde{\sigma}_{\alpha,\sigma,n}=T_{\alpha,\sigma}(P_{n})\text{ \ \ for
\ \ }T_{\alpha,\sigma}(Q)=\text{\textrm{argmin}}_{\tilde{\sigma}}%
M_{\alpha,\sigma}(Q,\tilde{\sigma}),\text{ }Q\in\mathcal{Q}%
\]
where%
\[
M_{\alpha,\sigma}(Q,\tilde{\sigma})=\tilde{M}_{\alpha,\sigma}(Q,\tilde{\sigma
}/\sigma)\text{ \ \ }%
\]
for (cf. (\ref{M}))%
\[
\tilde{M}_{\alpha,\sigma}(Q,s)=\frac{s^{\alpha}}{\left(  1-\alpha\right)
\sqrt{\alpha s^{2}+1-\alpha}}+\int\frac{s^{\alpha}}{\alpha}\exp\left\{
\frac{\alpha x^{2}\left[  s^{-2}-1\right]  }{2\sigma^{2}}\right\}
\mathrm{d}Q(x).
\]
$\medskip$

Put in accordance with (\ref{f00}) and (\ref{p3a})%
\begin{align}
\mathbf{\psi}_{\alpha,\sigma}(x,\tilde{\sigma})  & =\frac{\mathrm{d}%
}{\mathrm{d}\tilde{\sigma}}M_{\alpha,\sigma}(\delta_{x},\tilde{\sigma}%
)=\frac{1}{\sigma}\left(  \frac{\mathrm{d}}{\mathrm{d}s}\tilde{M}%
_{\alpha,\sigma}(\delta_{x},s)\right)  _{s=\tilde{\sigma}/\sigma}%
\medskip\nonumber\\
& =-\frac{1}{\sigma}\left[  s^{\alpha-1}\left(  \frac{\alpha\left(
s^{2}-1\right)  }{\left(  \alpha s^{2}+1-\alpha\right)  ^{3/2}}+\left[
\left(  \frac{x}{\sigma s}\right)  ^{2}-1\right]  e^{\alpha x^{2}\left[
s^{-2}-1\right]  /2\sigma^{2}}\right)  \right]  _{s=\tilde{\sigma}/\sigma
}\label{A5}\\
& =-\left(  \frac{\tilde{\sigma}}{\sigma}\right)  ^{\alpha-1}\left(
\frac{\alpha\left(  \tilde{\sigma}^{2}-\sigma^{2}\right)  }{\left[
\alpha\tilde{\sigma}^{2}+(1-\alpha)\sigma^{2}\right]  ^{3/2}}+\frac{1}{\sigma
}\left[  \left(  \frac{x}{\tilde{\sigma}}\right)  ^{2}-1\right]  e^{\alpha
x^{2}\left[  \tilde{\sigma}^{-2}-\sigma^{-2}\right]  /2}\right)  .\nonumber
\end{align}
By differentiating this expression with respect to $\tilde{\sigma}$\ and using
(\ref{f6}) we obtain the matrix%
\begin{equation}
I_{\alpha,\sigma}(\tilde{\sigma}):=\boldsymbol{I}_{\alpha,\sigma}%
(P_{\tilde{\sigma}})=\left(  \frac{\tilde{\sigma}}{\sigma}\right)  ^{\alpha
-1}\frac{2\sigma^{4}+\alpha^{2}(\tilde{\sigma}^{2}-\sigma^{2})^{2}}%
{\tilde{\sigma}[\alpha\tilde{\sigma}^{2}+(1-\alpha)\sigma^{2}]^{5/2}%
}.\label{a9}%
\end{equation}
Hence, by Theorem 2.1.1,\ the influence function of max\b{D}$_{\alpha}%
$-estimators at the data generating distributions $P_{\mathbb{\sigma}_{0}}$
are for all $0<\alpha<1$
\begin{align}
\text{IF}(x;\tilde{T}_{\alpha,\sigma},\mathbb{\sigma}_{0})  & =-\frac
{\mathbf{\psi}_{\alpha,\sigma}(x,\mathbb{\sigma}_{0})}{I_{\alpha,\sigma
}(\mathbb{\sigma}_{0})}\medskip\nonumber\\
& =\Delta_{\alpha,\sigma}(x;\mathbb{\sigma}_{0})+\frac{\alpha\sigma_{0}\left(
\mathbb{\sigma}_{0}^{2}-\sigma^{2}\right)  \left[  \alpha\mathbb{\sigma}%
_{0}^{2}+(1-\alpha)\sigma^{2}\right]  }{2\sigma^{4}+\alpha^{2}(\mathbb{\sigma
}_{0}^{2}-\sigma^{2})^{2}}\label{a10}%
\end{align}
where%
\begin{equation}
\Delta_{\alpha,\sigma}(x;\mathbb{\sigma}_{0})=\frac{\left[  \alpha
\mathbb{\sigma}_{0}^{2}+(1-\alpha)\sigma^{2}\right]  ^{5/2}\left[  \left(
x/\mathbb{\sigma}_{0}\right)  ^{2}-1\right]  \exp\left\{  \alpha x^{2}\left[
\mathbb{\sigma}_{0}^{-2}-\sigma^{-2}\right]  /2\right\}  }{\sigma\left[
2\sigma^{4}+\alpha^{2}(\mathbb{\sigma}_{0}^{2}-\sigma^{2})^{2}\right]
/\mathbb{\sigma}_{0}}.\label{a11}%
\end{equation}
This formula remains valid also for $\alpha=0$ since in this case (\ref{a10})
reduces to the well known influence function
\[
\text{IF}(x;MLE,\mathbb{\sigma}_{0})=\frac{\mathbb{\sigma}_{0}\left[  \left(
x/\mathbb{\sigma}_{0}\right)  ^{2}-1\right]  }{2}%
\]
obtained from the limit values%
\[
\mathbf{\psi}_{0,\sigma}(x,\mathbb{\sigma}_{0})=-\left[  \left(
x/\mathbb{\sigma}_{0}\right)  ^{2}-1\right]  /\mathbb{\sigma}_{0}\text{
\ \ and \ \ }I_{0,\sigma}(\tilde{\sigma})=2/\mathbb{\sigma}_{0}^{2}%
\]
which do not depend on the escort parameter . We see from the formula
(\ref{a11}) that the influence curve is unbounded for all $\sigma,\sigma
_{0}>0$ and $\alpha\geq0$. For $\alpha>0$\ and $\sigma\neq\sigma_{0}$\ we get
IF$(\mathbb{\sigma}_{0};\tilde{T}_{\alpha,\sigma},\mathbb{\sigma}_{0})\neq
0.$\ If moreover $\sigma<\sigma_{0}$\ then IF$(x;\tilde{T}_{\alpha,\sigma
},\mathbb{\sigma}_{0})$\ increases to infinity exponentially fast for
$|x|\rightarrow\infty.$ Thus $\tilde{T}_{\alpha,\sigma}$\ with $\alpha>0$\ and
$\sigma\neq\sigma_{0}$\ are strongly non-robust.

\paragraph{Example 2.1.4: Power subdivergence estimator in Pareto family.}

It is hard to find simpler nontrivial examples of the max\b{D}$_{\alpha}%
$-estimators\ than the estimators of location (\ref{63a}),\ (\ref{p9}) from
Example 2.1.2. Another relatively simple example is the family of max\b
{D}$_{\alpha}$-estimators\ in the Pareto model with the family of measures
${\mathcal{P}}=\{P_{\theta}:\theta>0\}$ defined on the interval ${\mathcal{X}%
}=(1,\infty)$ by the densities%
\begin{equation}
p_{\theta}(x)=\frac{\theta}{x^{\theta+1}}.\label{2.17}%
\end{equation}
with the mean values finite equal $\theta/(\theta-1)$\ in the domain
$\theta>1$\ and variances finite and equal $\theta/[(\theta-2)(\theta-1)^{2}%
]$\ in the domain $\theta>2.$ As before, the estimates $\tilde{\theta}%
_{\alpha,\theta,n}$\ depend on the divergence parameters $\alpha\geq0$\ and
escort parameters $\theta>0.$ By (\ref{MLE}), for $\alpha=0$ we get the MLE
estimates%
\[
\tilde{\theta}_{0,\theta,n}=\text{\textrm{argmax}}_{\tilde{\theta}}\text{
}\Sigma_{i=1}^{n}\ln p_{\tilde{\theta}}(X_{i})=\left(  \frac{1}{n}\sum
_{i=1}^{n}\ln X_{i}\right)  ^{-1}.
\]
For $0<\alpha<1$ we can use the criterion function%
\begin{equation}
M_{\alpha,\theta}(Q,\tilde{\theta})=\frac{1}{1-\alpha}P_{\tilde{\theta}}%
\cdot\left(  \frac{p_{\theta}}{p_{\tilde{\theta}}}\right)  ^{\alpha}+\frac
{1}{\alpha}Q\cdot\left(  \frac{p_{\theta}}{p_{\tilde{\theta}}}\right)
^{\alpha},\text{ \ \ }Q\in\mathcal{Q}\label{a6}%
\end{equation}
\ of (\ref{gg3}), or its derivative%
\begin{equation}
\Psi_{\alpha,\theta}(Q,\tilde{\theta})=\frac{\mathrm{d}}{\mathrm{d}%
\tilde{\theta}}M_{\alpha,\theta}(Q,\tilde{\theta})=P_{\tilde{\theta}}%
\cdot\left(  \frac{p_{\theta}}{p_{\tilde{\theta}}}\right)  ^{\alpha}%
s_{\tilde{\theta}}-Q\cdot\left(  \frac{p_{\theta}}{p_{\tilde{\theta}}}\right)
^{\alpha}s_{\tilde{\theta}}\label{a7}%
\end{equation}
given by (\ref{p3}), where in the present situation%
\[
P_{\tilde{\theta}}\cdot\left(  \frac{p_{\theta}(x)}{p_{\tilde{\theta}}%
(x)}\right)  ^{\alpha}=\frac{\theta^{\alpha}\tilde{\theta}^{1-\alpha}}%
{\alpha\theta+(1-\alpha)\tilde{\theta}},\ \ \text{and}\ \ \ s_{\theta
}(x)=\frac{1}{\theta}-\ln x.
\]
Substituting these expressions in (\ref{a6}), (\ref{a7}) we get the desired
asymptotic characteristics of the max\b{D}$_{\alpha}$-estimators\ $\tilde
{\theta}_{\alpha,\theta,n}$ obtained as argminima of the functions
$M_{\alpha,\theta}(P_{n},\tilde{\theta})$\ or, equivalently, as solutions of
the equations $\Psi_{\alpha,\theta}(P_{n},\tilde{\theta})=0$ in the variable
$\tilde{\theta}$. Further, by (\ref{f00}),%
\[
\mathbf{\psi}_{\alpha,\theta}(x,\tilde{\theta})=\Psi_{\alpha,\theta}%
(\delta_{x},\tilde{\theta})=P_{\tilde{\theta}}\cdot\left(  \frac{p_{\theta}%
}{p_{\tilde{\theta}}}\right)  ^{\alpha}s_{\tilde{\theta}}-\left(
\frac{p_{\theta}(x)}{p_{\tilde{\theta}}(x)}\right)  ^{\alpha}s_{\tilde{\theta
}}(x)
\]
and using Theorem 2.1.1 one easily obtains the influence functions of the
estimators $\tilde{\theta}_{\alpha,\theta,n}$ under consideration.

%%%%%

\subsection{Power superdivergence estimators}

In this subsection we deal with the $\min$\={D}$_{\alpha}$-estimators
$\theta_{\alpha,n}$ with the power parameters $\alpha\geq0.$\ For $\alpha
=0$\ they coincide with the MLE's (\ref{MLE}). Therefore we consider
$\alpha>0$ when these estimators are defined by (\ref{gg2}) and (\ref{gg3}).
Restrict ourselves for simplicity to $0<\alpha<1$\ and denote the function
$\Psi_{\alpha,\theta}(Q,\tilde{\theta})$ from (\ref{p3}) in previous
subsection temporarily by $\tilde{\Psi}_{\alpha,\theta}(Q,\tilde{\theta}%
),$\ i.e. let%
\[
\tilde{\Psi}_{\alpha,\theta}(Q,\tilde{\theta})=P_{\tilde{\theta}}\cdot\left(
\frac{p_{\theta}}{p_{\tilde{\theta}}}\right)  ^{\alpha}s_{\tilde{\theta}%
}-Q\cdot\left(  \frac{p_{\theta}}{p_{\tilde{\theta}}}\right)  ^{\alpha
}s_{\tilde{\theta}}.
\]
Further, let $\tilde{T}_{\alpha,\theta}(Q)$ be solution of the equation
$\tilde{\Psi}_{\alpha,\theta}(Q,\tilde{\theta})=0$\ in variable $\tilde
{\theta}$, i.e.
\begin{equation}
\tilde{\Psi}_{\alpha,\theta}(Q,\tilde{T}_{\alpha,\theta}(Q))=0\text{ \ \ for
all }\theta\in\Theta\text{.}\label{b3}%
\end{equation}
Finally, let $M_{\alpha,\theta}(Q,\tilde{T}_{\alpha,\theta}(Q))$ be the
function of variable $\theta\in\Theta$ obtained by inserting $\tilde{\theta
}=\tilde{T}_{\alpha,\theta}(Q)$\ in the function $M_{\alpha,\theta}%
(Q,\tilde{\theta})$\ defined in (\ref{gg3}). According to (\ref{gg2}) and
(\ref{gg3}), the maximizers%
\begin{equation}
T_{\alpha}(Q)=\text{\textrm{argmax}}_{\theta}\text{ }M_{\alpha,\theta
}(Q,\tilde{T}_{\alpha,\theta}(Q))\label{b0}%
\end{equation}
generate the $\min$\={D}$_{\alpha}$-estimators $\theta_{\alpha,n}$ under
consideration in the sense that $\theta_{\alpha,n}=T_{\alpha}(P_{n}).\medskip$

In the following theorem we consider the score function $s_{\theta}%
=\mathring{p}_{\theta}/p_{\theta}$\ and we put for brevity $\tilde{\tau
}_{\alpha,\theta}=\tilde{T}_{\alpha,\theta}(Q)$.

\paragraph{Theorem 2.2.1.}

For all $0<\alpha<1$\ the maximizers (\ref{b0}) solve the equations
$\Psi_{\alpha}(Q,\theta)=0$\ in variable $\theta\in\Theta$\ for the function%
\begin{equation}
\Psi_{\alpha}(Q,\theta)=\frac{\mathrm{d}}{\mathrm{d}\theta}\text{ }%
M_{\alpha,\theta}(Q,\tilde{\tau}_{\alpha,\theta})=\frac{\alpha}{1-\alpha
}P_{\tilde{\tau}_{\alpha,\theta}}\cdot\left(  \frac{p_{\theta}}{p_{\tilde
{\tau}_{\alpha,\theta}}}\right)  ^{\alpha}s_{\theta}+Q\cdot\left(
\frac{p_{\theta}}{p_{\tilde{\tau}_{\alpha,\theta}}}\right)  ^{\alpha}%
s_{\theta}.\medskip\label{b6}%
\end{equation}
Consequently the corresponding $\min$\={D}$_{\alpha}$-estimators
$\theta_{\alpha,n}=T_{\alpha}(P_{n})$\ are solutions of the equations
\begin{equation}
\frac{\alpha}{1-\alpha}P_{\tilde{\tau}_{\alpha,\theta}}\cdot\left(
\frac{p_{\theta}}{p_{\tilde{\tau}_{\alpha,\theta}}}\right)  ^{\alpha}%
s_{\theta}+\frac{1}{n}\sum_{i=1}^{n}\left(  \frac{p_{\theta}(X_{i})}%
{p_{\tilde{\tau}_{\alpha,\theta}}(X_{i})}\right)  ^{\alpha}s_{\theta}%
(X_{i})=0.\label{b2}%
\end{equation}
$\medskip$

\paragraph{Proof.}

By (\ref{gg3})%
\begin{equation}
M_{\alpha,\theta}(Q,\tilde{\theta})=\frac{1}{1-\alpha}P_{\tilde{\theta}}%
\cdot\left(  \frac{p_{\theta}}{p_{\tilde{\theta}}}\right)  ^{\alpha}+\frac
{1}{\alpha}Q\cdot\left(  \frac{p_{\theta}}{p_{\tilde{\theta}}}\right)
^{\alpha}\nonumber
\end{equation}
so that%
\begin{align*}
\frac{\mathrm{d}}{\mathrm{d}\theta}\text{ }M_{\alpha,\theta}(Q,\tilde{\tau
}_{\alpha,\theta})  & =\left(  \frac{\mathrm{d}}{\mathrm{d}\theta}\text{
}M_{\alpha,\theta}(Q,\tilde{\theta})\right)  _{\tilde{\theta}=\tilde{\tau
}_{\alpha,\theta}}+\left(  \frac{\mathrm{d}}{\mathrm{d}\theta}\text{
}M_{\alpha,\tilde{\theta}}(Q,\tilde{\tau}_{\alpha,\theta})\right)
_{\tilde{\theta}=\theta}\text{ }\medskip\\
& =\frac{\alpha}{1-\alpha}P_{\tilde{\tau}_{\alpha,\theta}}\cdot\left(
\frac{p_{\theta}}{p_{\tilde{\tau}_{\alpha,\theta}}}\right)  ^{\alpha}%
s_{\theta}+Q\cdot\left(  \frac{p_{\theta}}{p_{\tilde{\tau}_{\alpha,\theta}}%
}\right)  ^{\alpha}s_{\theta}\\
& +\left(  \frac{\mathrm{d}}{\mathrm{d}\tau}M_{\alpha,\theta}(Q,\tau)\right)
_{\tau=\tilde{\tau}_{\alpha,\theta}}\text{.}\frac{\mathrm{d}\tilde{\tau
}_{\alpha,\theta}}{\mathrm{d}\theta}\medskip\\
& =\frac{\alpha}{1-\alpha}P_{\tilde{\tau}_{\alpha,\theta}}\cdot\left(
\frac{p_{\theta}}{p_{\tilde{\tau}_{\alpha,\theta}}}\right)  ^{\alpha}%
s_{\theta}+Q\cdot\left(  \frac{p_{\theta}}{p_{\tilde{\tau}_{\alpha,\theta}}%
}\right)  ^{\alpha}s_{\theta}+\tilde{\Psi}_{\alpha,\theta}(Q,\tilde{\tau
}_{\alpha,\theta})\text{.}\frac{\mathrm{d}\tilde{\tau}_{\alpha,\theta}%
}{\mathrm{d}\theta}.
\end{align*}
Using (\ref{b3}) we obtain (\ref{b6}) and (\ref{b2}).\hfill$\blacksquare
\medskip$

\paragraph{Corollary 2.2.1.}

The influence functions IF$(x;T_{\alpha},\theta)$\ of all $\min$\={D}%
$_{\alpha}$-estimators $\theta_{\alpha,n}=T_{\alpha}(P_{n})$ with power
parameters $0<\alpha<1$\ at $P_{\theta}\in\mathcal{P}$ coincide with the
influence function
\begin{equation}
\text{IF}(x;T_{0},\theta)=\boldsymbol{I}(\theta)^{-1}s_{\theta}(x)\text{
\ \ \ \ \ \ \ \ \ \ \ \ \ \ \ \ \ \ \ \ (cf. (\ref{f66}) and (\ref{if}%
))}\label{b5}%
\end{equation}
of the MLE $\theta_{0,n}=T_{0}(P_{n})$.

\paragraph{Proof.}

By Theorem 2.2, the max\b{D}$_{\alpha}$-estimators $\tilde{\theta}%
_{\alpha,\theta n}=\tilde{T}_{\alpha,\theta}(P_{n})$ are Fisher consistent.
Hence for $Q=P_{\theta_{0}}\ $we get $\tilde{\tau}_{\alpha,\theta}:=\tilde
{T}_{\alpha,\theta}(P_{\theta_{0}})=\theta_{0}$\ in (\ref{b6}). Consequently
it follows from (\ref{f00}) and (\ref{b6}) that the $\boldsymbol{\psi}%
$-functions%
\[
\boldsymbol{\psi}_{\alpha}(x,\tilde{\tau}_{\alpha,\theta})\equiv\Psi_{\alpha
}(\delta_{x},\tilde{\tau}_{\alpha,\theta})=\frac{\alpha}{1-\alpha}%
P_{\tilde{\tau}_{\alpha,\theta}}\cdot\left(  \frac{p_{\theta_{0}}}%
{p_{\tilde{\tau}_{\alpha,\theta}}}\right)  ^{\alpha}s_{\theta_{0}}+\delta
_{x}\cdot\left(  \frac{p_{\theta_{0}}}{p_{\tilde{\tau}_{\alpha,\theta}}%
}\right)  ^{\alpha}s_{\theta_{0}}%
\]
of these estimators reduce for all $0<\alpha<1$\ to the score function
$s_{\theta_{0}}(x)$\ which is the $\boldsymbol{\psi}$-function of MLE $T_{0}
$. Similarly, we get from (\ref{f66}) and (\ref{f6}) for all $0<\alpha<1$\ the
matrix $\boldsymbol{I}(\theta_{0})=P_{\theta_{0}}\cdot s_{\theta_{0}%
}^{\text{t}}s_{\theta_{0}}$ corresponding to the MLE. Therefore the influence
functions of all $\min$\={D}$_{\alpha}$-estimators under
considerations\ reduce to the influence MLE function (\ref{b5})\ which
completes the proof.\hfill$\blacksquare\medskip$

Formulas for the $\min$\={D}$_{\alpha}$-estimators of the normal location
and/or scale are seen from the examples of Subsection 2.1.

%%%%%%%%%%%%%%%%%%%%%%%%%%%%%%%%%%%%%%%%%%%%%%%%%

\section{DECOMPOSABLE PSEUDODISTANCES}

The $\phi$-divergences $D_{\phi}(P,Q)$, $\phi\in\mbox{\boldmath$\Phi$}$ can be
characterized by the \textsl{information processing property}, i.\thinspace e.
by the complete invariance w.r.t. the statistically sufficient transformations
of the observation space $({\mathcal{X}},{\mathcal{A}})$. This property is
useful but probably not unavoidable in the minimum distance estimation based
on similarity between theoretical and empirical distributions. Hence we admit
in the rest of the paper general \textsl{pseudodistances} $\mathfrak{D}(P,Q)$
which may not satisfy the information processing property.

\paragraph{Definition 3.1.}

We say that $\mathfrak{D}:$\ $\mathcal{P}\otimes\mathcal{P}^{+}\mapsto
\mathbb{R}$ is a pseudodistance of probability measures $P\in\mathcal{P}%
=\{P_{\theta}:\theta\in\Theta\}$ and $Q\in\mathcal{P}^{+}$if%
\begin{equation}
\mathfrak{D}(P_{\theta},P_{{\tilde{\theta}}})\geq0\text{ \ for all }%
\theta,{\tilde{\theta}}\in\Theta\text{\ with \ }\mathfrak{D}(P_{\theta
},P_{{\tilde{\theta}}})=0\text{ \ \ iff \ }\theta={\tilde{\theta}}.\label{24a}%
\end{equation}

\medskip

An additional restriction imposed in this section on pseudodistances
$\mathfrak{D}(P,Q)$ will be the \textsl{decomposability}.

\paragraph{Definition 3.2.}

A pseudodistance $\mathfrak{D}$ on $\mathcal{P}\otimes\mathcal{P}^{+}$is a
$\boldsymbol{decomposable}$ if there exist functionals $\mathfrak{D}%
^{0}:\mathcal{P}\mapsto\mathbb{R}$, $\mathfrak{D}^{1}:\mathcal{P}^{+}%
\mapsto\mathbb{R}$ and measurable mappings%
\begin{equation}
\rho_{\theta}:\mathcal{X\mapsto}\mathbb{R},\text{ \ \ \ }\theta\in
\Theta\label{24aa}%
\end{equation}
\ such that for all $\theta\in\Theta$\ and $Q\in\mathcal{P}^{+}$ the
expectations $Q\cdot\rho_{\theta}$\ exist and%
\begin{equation}
\mathfrak{D}(P_{\theta},Q)=\mathfrak{D}^{0}(P_{\theta})+\mathfrak{D}%
^{1}(Q)+Q\cdot\rho_{\theta}.\label{26a}%
\end{equation}

\paragraph{Definition 3.3.}

We say that a functional $T_{\mathfrak{D}}:\mathcal{Q}\mapsto\Theta$\ for
$\mathcal{Q}=\mathcal{P}^{+}\cup\mathcal{P}_{\text{emp}}$\ defines a
$\boldsymbol{minimum}$ $\boldsymbol{pseudodistance}$ $\boldsymbol{estimator}$
(briefly, $\min\mathfrak{D}${-estimator})if $\mathfrak{D}(P_{\theta},Q)$\ is a
decomposable pseudodistance on $\mathcal{P}\otimes\mathcal{P}^{+}$\ and\ the
parameters $T_{\mathfrak{D}}(Q)\in\Theta$\ minimize $\mathfrak{D}%
^{0}(P_{\theta})+Q\cdot\rho_{\theta}$ on $\Theta,$\ in symbols
\begin{equation}
T_{\mathfrak{D}}(Q)=\mathrm{argmin}_{\theta}\left[  \mathfrak{D}^{0}%
(P_{\theta})+Q\cdot\rho_{\theta}\right]  \text{ \ \ \ for all \ }%
Q\in\mathcal{Q}.\label{c11}%
\end{equation}
In particular, for $Q=P_{n}\in\mathcal{P}_{\text{emp}}$\ \
\begin{equation}
\theta_{\mathfrak{D},n}:=T_{\mathfrak{D}}(P_{n})=\mathrm{argmin}_{\theta
}\left[  \mathfrak{D}^{0}(P_{\theta})+\frac{1}{n}\sum_{i=1}^{n}\rho_{\theta
}(X_{i})\right]  \text{ \ \ if \ }P_{n}=\frac{1}{n}\sum_{i=1}^{n}\delta
_{X_{i}}.\label{c1}%
\end{equation}
\medskip

\paragraph{Theorem 3.1.}

Every $\min\mathfrak{D}$-estimator\
\begin{equation}
\theta_{\mathfrak{D},n}=\mathrm{argmin}_{\theta}\left[  \mathfrak{D}%
^{0}(P_{\theta})+\frac{1}{n}\sum_{i=1}^{n}\rho_{\theta}(X_{i})\right]
\label{c2}%
\end{equation}
is Fisher consistent in the sense that%
\begin{equation}
T_{\mathfrak{D}}(P_{\theta_{0}})=\mathrm{argmin}_{\theta}\mathfrak{D}%
(P_{\theta},P_{\theta_{0}})=\theta_{0}\text{ \ \ for all }\theta_{0}\in
\Theta.\label{c3}%
\end{equation}

\paragraph{Proof.}

Consider arbitrary fixed $\theta_{0}\in\Theta$. Then, by assumptions,
$\mathfrak{D}^{1}(P_{\theta_{0}})$ is a finite constant. Therefore (\ref{c11})
together the definition of pseudodistance implies%
\begin{align*}
T_{\mathfrak{D}}(P_{\theta_{0}})  & =\mathrm{argmin}_{\theta}\left[
\mathfrak{D}^{0}(P_{\theta})+Q\cdot\rho_{\theta}\right] \\
& =\mathrm{argmin}_{\theta}\left[  \mathfrak{D}^{0}(P_{\theta})+\mathfrak{D}%
^{1}(P_{\theta_{0}})+Q\cdot\rho_{\theta}\right] \\
& =\mathrm{argmin}_{\theta}\mathfrak{D}(P_{\theta},P_{\theta_{0}})=\theta_{0}.
\end{align*}
\hfill$\blacksquare$\bigskip

The decomposability of pseudodistance $\mathfrak{D}(P_{\theta},Q)$ leads to
the additive structure of the criterion
\begin{equation}
\mathfrak{D}(P_{\theta},P_{n})\sim\mathfrak{D}^{0}(P_{\theta})+P_{n}\cdot
\rho_{\theta}=\mathfrak{D}^{0}(P_{\theta})+\frac{1}{n}\sum_{i=1}^{n}%
\rho_{\theta}(X_{i})\label{c4}%
\end{equation}
in the definition (\ref{c2}) of the $\min\mathfrak{D}$-estimators which opens
the possibility to apply the methods of the asymptotic theory of
$M$-estimators (cf. Hampel et al. (1986), van der Vaart and Wellner (1996),
van der Vaart (1998) or Mieske and Liese (2008)). $\medskip$

The general $\min\mathfrak{D}$-estimators and their special classes studied in
Subsections 3.1, 3.2 below were introduced in Vajda (2008). They contain as a
subclass all the $\max$\b{D}$_{\phi}$-estimators of Section~2. To see this
suppose that the assumptions of Section~2 related to the estimators (\ref{c5})
hold and consider for arbitrary fixed $(\phi,\tau)\in
\mbox{\boldmath$\Phi$}\otimes\Theta$ the well defined expressions%
\[
\mathfrak{D}_{\phi,\tau}^{0}(P_{\theta})=-\text{ }P_{\tau}\cdot\phi^{\prime
}\left(  {\frac{p_{\tau}}{p_{\theta}}}\right)  ,\text{ \ \ \ }\rho_{\phi
,\tau,\theta}=-\text{ }\phi^{\sharp}\left(  {\frac{p_{\tau}}{p_{\theta}}%
}\right)  \medskip
\]
and
\[
\mathfrak{D}_{\phi,\tau}^{1}(Q)=-\text{ }\inf_{\theta}\left[  \mathfrak{D}%
_{\phi,\tau}^{0}(P_{\theta})+Q\cdot\rho_{\phi,\tau,\theta}\right]  .
\]

\paragraph{Theorem 3.2.}

The sum
\begin{equation}
\mathfrak{D}(P_{\theta},Q):=\mathfrak{D}_{\phi,\tau}^{0}(P_{\theta
})+\mathfrak{D}_{\phi,\tau}^{1}(Q)+Q\cdot\rho_{\phi,\tau,\theta}\label{c6}%
\end{equation}
is a pseudodistance on $\mathcal{P}\otimes${$\mathcal{P}^{+}$} and the maximum
subdivergence estimator
\begin{equation}
\theta_{\phi,\tau,n}=\mbox{argmax}_{\theta}\left[  P_{\tau}\cdot\phi^{\prime
}\left(  {\frac{p_{\tau}}{p_{\theta}}}\right)  +{\frac{1}{n}}\sum_{i=1}%
^{n}\phi^{\sharp}\left(  {\frac{p_{\tau}(X_{i})}{p_{\theta}(X_{i})}}\right)
\right] \label{c5}%
\end{equation}
of Section~2 with the divergence parameter $\phi\in\mbox{\boldmath$\Phi$}$ and
escort parameter $\tau\in\Theta$ is the $\min\mathfrak{D}$-estimator for the
decomposable pseudodistance (\ref{c6}).

\paragraph{Proof.}

Fix $(\phi,\tau)\in\mbox{\boldmath$\Phi$}\otimes\Theta$ and let the
assumptions of Section~2 related to the estimators (\ref{c5}) hold. Then for
any $\theta_{0}\in\Theta$
\[
\mathfrak{D}(P_{\theta_{0}},Q)=\mathfrak{D}_{\phi,\tau}^{0}(P_{\theta_{0}%
})+Q\cdot\rho_{\phi,\tau,\theta_{0}}-\inf_{\theta}\left[  \mathfrak{D}%
_{\phi,\tau}^{0}(P_{\theta_{0}})+Q\cdot\rho_{\phi,\tau,\theta_{0}}\right]
\geq0.
\]
If $Q\in{\mathcal{P}}$ then, by (\ref{17a}) and (\ref{17g}),
\begin{align*}
\mathfrak{D}_{\phi,\tau}(P_{\theta_{0}},Q)  & =\sup_{\theta}\left[
P_{\theta_{0}}\cdot\phi^{\prime}\left(  {\frac{p_{\tau}}{p_{\theta_{0}}}%
}\right)  +Q\cdot\phi^{\sharp}\left(  {\frac{p_{\tau}}{p_{\theta}}}\right)
\right]  -P_{\theta_{0}}\cdot\phi^{\prime}\left(  {\frac{p_{\tau}}%
{p_{\theta_{0}}}}\right)  +Q\cdot\phi^{\sharp}\left(  {\frac{p_{\tau}%
}{p_{\theta_{0}}}}\right) \\
& \\
& =D_{\phi}(P_{\theta_{0}},Q)-\text{\b{D}}_{\phi,\tau}(P_{\theta_{0}},Q).
\end{align*}
By Theorem~2.1, this difference is zero if and only if $Q=P_{\theta_{0}}$
which proves that (\ref{c6}) is pseudodistance on $\mathcal{P}\otimes
${$\mathcal{P}^{+}$}$.$ On the other hand, obviously, (\ref{c5}) satisfies%
\[
\theta_{\phi,\tau,n}=\mbox{argmin}_{\theta}\left[  \mathfrak{D}_{\phi,\tau
}^{0}(P_{\theta})+P_{n}\cdot\rho_{\phi,\tau,\theta}\right]
\]
so that it is $\min\mathfrak{D}$-estimator for the pseudodistance (\ref{c6})
which completes the proof.\hfill$\blacksquare$\bigskip

The minimum superdivergence estimators $\theta_{\phi,n}$ of Section~2 (the
$\min$\={D}$_{\phi}$-estimators) minimize the suprema
\[
\sup_{\tau}\mathfrak{D}(P_{\theta},Q)\ \ \ \ \ \mbox{for}\ Q=P_{n}%
\]
of the decomposable pseudodistance (\ref{c6}). However, the suprema of
decomposable pseudodistances are not in general decomposable pseudodistances.
Therefore the standard theory of $M$-estimators is not applicable to this
class of estimators. An exception is the MLE $\theta_{\phi_{0},n}$ obtained
for the logarithmic function $\phi_{0}$ given in~(\ref{7}).

%%%%%%%%%%%%%%%%%%%%

\subsection{Power pseudodistance estimators}

In this subsection we study a special class of pseudodistances $\mathfrak{D}%
_{\psi}(P_{\theta},Q)$\ defined on $\mathcal{P}\otimes\mathcal{P}^{+}$\ by the
integral formula
\begin{equation}
\mathfrak{D}_{\psi}(P_{\theta},Q)=\int\psi(p_{\theta},q)\,\mathrm{d}%
\lambda\text{ \ \ \ for\ }p_{\theta}=\frac{\mathrm{d}P_{\theta}}%
{\mathrm{d}\lambda},q=\frac{\mathrm{d}Q}{\mathrm{d}\lambda}\label{25}%
\end{equation}
where $\psi(s,t)$ are \textsl{reflexive} in the sense that they are
nonnegative functions of arguments $s,t>0$ with $\psi(s,t)=0$ iff $s=t$. If a
function $\psi$\ is reflexive and also \textsl{decomposable} in the sense%
\begin{equation}
\psi(s,t)=\psi^{0}(s)+\psi^{1}(t)+\rho(s)\,t,\text{ \ \ }s,t\geq0\label{27}%
\end{equation}
for some $\psi^{0},\psi^{1},\rho:(0,\infty)\rightarrow\mathbb{R}$ then the
corresponding $\psi$-pseudodistance (\ref{25}) is a \textsl{decomposable
pseudodistance} satisfying%
\begin{equation}
\mathfrak{D}_{\psi}(P_{\theta},Q)=\mathfrak{D}_{\psi}^{0}(P_{\theta
})+\mathfrak{D}_{\psi}^{1}(Q)+Q\cdot\rho_{\theta}\text{ \ \ \ (cf.
(\ref{26a}))}\label{27a}%
\end{equation}
for%
\begin{equation}
\mathfrak{D}_{\psi}^{0}(P_{\theta})=\int\psi^{0}(p_{\theta})\,\mathrm{d}%
\lambda,\text{ \ }\mathcal{D}_{\psi}^{1}(Q)=\int\psi^{1}(q)\,\mathrm{d}%
\lambda\text{\ \ and \ }\rho_{\theta}=\rho(p_{\theta}).\label{27b}%
\end{equation}

\paragraph{Example 3.1.1.}

The $\phi$-divergences $D_{\phi}(P_{\theta},Q)$ are special $\psi
$-pseudodistances (\ref{25}) for the functions%
\begin{equation}
\psi(s,t)=\phi(s/t)\,t-\phi^{\prime}(1)(s-t),\quad s,t>0\label{26}%
\end{equation}
since they are nonnegative and reflexive, and (\ref{26}) implies
$\mathfrak{D}_{\psi}(P_{\theta},Q)=D_{\phi}(P_{\theta},Q)$ for all
$P\in\mathcal{P},Q\in\mathcal{P}^{+}$ when $\phi\in\mathbf{\Phi}$\ and $\psi
$\ are related by (\ref{26}). However, the functions (\ref{26}) in general do
not satisfy the decomposability condition (\ref{27}) so that the $\phi
$-divergences are\ not in general decomposable pseudodistances. An exception
is the logarithmic function $\phi=\phi_{0}$\ defined in (\ref{7}) for which
the $\min\mathfrak{D}_{\phi_{0}}$-estimator is the MLE.

\paragraph{Example 3.1.2: $L_{2}$-estimator}

The quadratic function $\psi(s,t)=(s-t)^{2}\ $is reflexive and also
decomposable in the sense of (\ref{27}). Thus it defines the decomposable
pseudodistance
\[
\mathfrak{D}_{\psi}(P_{\theta},Q)=\int(p_{\theta}-q)^{2}\,\mathrm{d}%
\lambda=\Vert p_{\theta}-q\Vert^{2}%
\]
on $\mathcal{P}\otimes\mathcal{P}^{+}$ for $\mathcal{P}^{+}\subset
L_{2}(\lambda)$. It is easy to verify that the decomposability in the sense of
(\ref{27a}) holds for%
\[
\mathfrak{D}_{\psi}^{0}(P_{\theta})=\int p_{\theta}^{2}\,\mathrm{d}%
\lambda\text{\ \ \ }\mathcal{D}_{\psi}^{1}(Q)=\int q^{2}\,\mathrm{d}%
\lambda_{Q},\text{ \ and \ \ }\rho_{\theta}=-2p_{\theta}.
\]
The corresponding $\min\mathfrak{D}_{\psi}$-estimator defined by (\ref{c2}) is
in this case the $L_{2}$-estimator
\begin{equation}
{\theta}_{n}=\mbox{argmin}_{{\theta}}\left[  \int p_{{\theta}}^{2}%
\,\mathrm{d}\lambda-{\frac{2}{n}}\sum_{i=1}^{n}p_{{\theta}}(X_{i})\right]
\label{29}%
\end{equation}
which is known to be robust but not efficient (see e.g. Hampel et al.
(1986)).\bigskip

To build a smooth bridge between the robustness and efficiency, one needs to
replace the reflexive and decomposable functions $\psi$ by families
$\{\psi_{\alpha}:\alpha\geq0\}$ of reflexive functions decomposable in the
sense
\begin{equation}
\psi_{\alpha}(s,t)=\psi_{\alpha}^{0}(s)+\psi_{\alpha}^{1}(t)+\rho_{\alpha
}(s)\,t\quad\text{for all }\alpha\geq0\text{\ \ \ }\mbox{(cf.
(\ref{27}))}\label{30}%
\end{equation}
with the limits at satisfying for some constant $\varkappa$\ all $s>0$\ the
conditions%
\begin{equation}
\psi_{0}^{0}(s)=\lim_{\alpha\downarrow0}\psi_{\alpha}^{0}(s)=\varkappa\text{
}s\text{\ \ and \ }\lim_{\alpha\downarrow0}\rho_{\alpha}(s)=\rho_{0}(s)=-\ln
s.\label{31}%
\end{equation}
Then for all $\alpha\geq0$ and $(P_{\theta},Q)\in\mathcal{P}\otimes
\mathcal{P}^{+}$ the family of $\psi_{\alpha}$-pseudodistances%
\begin{equation}
\mathfrak{D}_{\alpha}(P_{\theta},Q):=\mathfrak{D}_{\psi_{\alpha}}(P_{\theta
},Q),\text{ \ \ \ }\alpha\geq0\label{31aa}%
\end{equation}
satisfies the decomposability condition\
\begin{equation}
\mathfrak{D}_{\alpha}(P_{\theta},Q)=\mathfrak{D}_{\alpha}^{0}(Q)+\mathfrak{D}%
_{\alpha}^{1}(P_{\theta})+Q\cdot\rho_{\alpha,\theta}\text{ \ \ \ (cf.
(\ref{26a}))}\label{31a}%
\end{equation}
for%
\begin{equation}
\mathfrak{D}_{\alpha}^{0}(P_{\theta})=\int\psi_{\alpha}^{0}(p_{\theta
})\,\mathrm{d}\lambda,\text{ \ \ }\mathfrak{D}_{\alpha}^{1}(Q)=\int
\psi_{\alpha}^{1}(q)\,\mathrm{d}\lambda\text{\ \ \ and \ }\rho_{\alpha,\theta
}=\rho_{\alpha}(p_{\theta}).\label{31b}%
\end{equation}
In other words, the pseudodistances $\mathfrak{D}_{\alpha}(P_{\theta},Q)$
defined by (\ref{31aa}) are decomposable and define in accordance with
(\ref{c2}) the family of $\min\mathfrak{D}_{\alpha}$-estimators
\begin{align}
{\theta}_{\alpha,n}  & =\mathrm{arg\min}_{{\theta}}\left[  \mathfrak{D}%
_{\psi_{\alpha}}^{0}(P_{{\theta}})+P_{n}\cdot\rho_{\alpha,\theta}\right]
\medskip\label{31c}\\
& =\mathrm{argmin}_{{\theta}}\left[  \int\psi_{\alpha}^{0}(p_{{\theta}%
})\,\mathrm{d}\lambda+\frac{1}{n}\sum_{i=1}^{n}\rho_{\alpha}(p_{{\theta}%
}(X_{i}))\right]  ,\text{ \ \ }\alpha\geq0.\label{31d}%
\end{align}
Here (\ref{31}) guarantees that this family contains as a special case for
$\alpha=0$ the efficient but non-robust MLE%
\begin{equation}
{\theta}_{0,n}=\mbox{argmin}_{{\theta}}\left[  \mathrm{const}-{\frac{1}{n}%
}\sum_{i=1}^{n}\ln p_{{\theta}}(X_{i})\right] \label{33}%
\end{equation}
while for $\alpha>0$\ the ${\theta}_{\alpha,n}$'s\ are expected to be less
efficient but more robust than ${\theta}_{0,n}$.$\medskip$

The rest of this subsection studies special family of decomposable
pseudodistances $\mathfrak{D}_{\alpha}(P_{\theta},Q).$\textsl{\ }It
is\textsl{\ }defined on $\mathcal{P}\otimes\mathcal{Q}$\ in accordance with
(\ref{31aa}) and (\ref{25}) by the functions
\begin{equation}
\psi_{\alpha}(s,t)=t^{1+\alpha}\left[  \alpha\phi_{1+\alpha}\left(  \frac
{s}{t}\right)  +(1-\alpha)\phi_{\alpha}\left(  \frac{s}{t}\right)  \right]
,\text{ \ \ }\alpha\geq0\label{34}%
\end{equation}
of variables $s,\,t>0$ where $\phi_{1+\alpha}$\ and $\phi_{\alpha}$\ are the
power functions defined by (\ref{6}), (\ref{7}). These functions satisfy
(\ref{30}), (\ref{31}) as it is clarified by the next theorem. In this theorem
and in the sequel we use for the function (\ref{34}) the relations%
\begin{align}
\psi_{\alpha}(s,t)  & =\frac{s^{1+\alpha}}{1+\alpha}+t^{1+\alpha}\left(
\frac{1}{\alpha}-\frac{1}{1+\alpha}\right)  -\frac{ts^{\alpha}}{\alpha
}\label{344}\\
& \nonumber\\
& =\frac{s^{1+\alpha}-t^{1+\alpha}}{1+\alpha}+t\left(  \frac{t^{\alpha}%
-1}{\alpha}-\frac{s^{\alpha}-1}{\alpha}\right) \label{345}%
\end{align}
when $\alpha>0$\ and%
\begin{align}
\psi_{0}(s,t)  & =s-t+t\ln t-t\ln s\label{346}\\
& \nonumber\\
& =\lim_{\alpha\downarrow0}\frac{s^{1+\alpha}-t^{1+\alpha}}{1+\alpha}+t\left(
\frac{t^{\alpha}-1}{\alpha}-\frac{s^{\alpha}-1}{\alpha}\right) \label{347}%
\end{align}
when $\alpha=0$.

%%%%%

\paragraph{Theorem 3.1.1.}

The power functions (\ref{34}) are reflexive and decomposable in the sense of
(\ref{30}) with
\begin{equation}
\psi_{\alpha}^{0}(s)=\frac{s^{1+\alpha}}{1+\alpha},\text{ \ }\psi_{\alpha}%
^{1}(t)=\left\{
\begin{array}
[c]{ll}%
t\left[  \frac{{\Large t}^{{\Large \alpha}}-{\large 1}}{{\Large \alpha}}%
-\frac{{\Large t}^{{\Large \alpha}}}{1+{\Large \alpha}}\right]  & \medskip\\
t\ln t-t &
\end{array}
\right.  \text{ and \ }\rho_{\alpha}(s)=\left\{
\begin{array}
[c]{ll}%
-{\frac{{\Large s}^{{\Large \alpha}}-1}{{\Large \alpha}\,}} &
\mbox{if}\ \alpha>0\medskip\\
-\ln s & \mbox{if}\text{ }\alpha=0.
\end{array}
\right. \label{34A}%
\end{equation}
Moreover, this family is continuous in the parameter $\alpha\downarrow0$ and
satisfies (\ref{31}) for $\varkappa=1$.

\paragraph{Proof.}

Decomposition (\ref{30}) for function $\psi_{\alpha}(s,t)$\ of (\ref{34}) into
the components (\ref{34A}) is clear from (\ref{345}) when $\alpha>0$\ and
(\ref{346}) when $\alpha=0$. The continuity in the parameter $\alpha
\downarrow0$ and (\ref{31}) for $\varkappa=1$ follow from (\ref{347}). We
shall prove the nonnegativity and reflexivity. For arbitrary arguments
$s,\,t>0$ and fixed parameters $a,\,b>0$ with the property $1/a+1/b=1$ it
holds
\begin{equation}
st\leq{\frac{s^{a}}{a}}+{\frac{t^{b}}{b}}\label{34a}%
\end{equation}
where = takes place iff $s^{a}=t^{b}$. Indeed, from the strict concavity of
the logarithmic function we deduce the inequality
\[
\ln(st)={\frac{1}{a}}\ln s^{a}+{\frac{1}{b}}\ln t^{b}\leq\ln\left(
{\frac{s^{a}}{a}}+{\frac{t^{b}}{b}}\right)
\]
and the stated condition for equality. Substituting $s\rightarrow s^{\alpha},$
$a\rightarrow(1+\alpha)/\alpha$ and $b\rightarrow1+\alpha$\ for $\alpha>0$ we
get%
\[
s^{\alpha}t\leq{\frac{s^{1+\alpha}}{(1+\alpha)/\alpha}}+{\frac{t^{1+\alpha}%
}{1+\alpha}}%
\]
with the equality condition $s^{\alpha a}=t^{b},$ i.e. $s^{1+\alpha
}=t^{1+\alpha}$.\ This implies that\ the function $\psi_{\alpha}(s,t)$\ is
nonnegative and reflexive. $\blacksquare$\bigskip

By (\ref{31aa}), (\ref{25}) and Theorem 3.1.1, the power functions (\ref{34})
generate%
\begin{equation}
\psi^{0}(p_{\theta})=\frac{{\small 1}}{{\small 1+\alpha}}p_{\theta}^{\alpha
}\text{ \ \ and \ \ }\rho_{\alpha}(p_{\theta})=\left\{
\begin{array}
[c]{ll}%
-{\frac{{\Large 1}}{{\Large \alpha}}}p_{\theta}^{\alpha} & \mbox{if}\ \alpha
>0\medskip\\
-\ln p_{\theta} & \mbox{if}\text{ }\alpha=0.
\end{array}
\right. \label{28a}%
\end{equation}
and define the family of decomposable pseudodistances%
\begin{equation}%
\begin{array}
[c]{ll}%
\mathfrak{D}_{\alpha}(P_{\theta},Q) & =\int\psi_{\alpha}(p_{\theta
},q)\,\mathrm{d}\lambda\medskip\\
& =\left\{
\begin{array}
[c]{ll}%
\frac{{\Large 1}}{{\Large 1+\alpha}}P_{\theta}\cdot p_{\theta}^{\alpha}%
+{\frac{{\Large 1}}{{\Large \alpha(1+\alpha})}}Q\cdot q^{\alpha}%
-{\frac{{\Large 1}}{{\Large \alpha}}Q\cdot}p_{\theta}^{\alpha} &
\mbox{if}\ \alpha>0\medskip\\
Q\cdot\left(  \ln q-\ln p_{\theta}\right)  & \mbox{if}\text{ }\alpha=0
\end{array}
\right.  \medskip
\end{array}
\label{36aa}%
\end{equation}
in (\ref{31d}). Relation of this family to the family of power divergences
$D_{\alpha}(P_{\theta},Q)$\ defined by (\ref{5}) is rigorously established in
the next theorem. It refers to the auxiliary family of functions%
\begin{equation}
\varphi_{\alpha}(s,t)=t\left[  \alpha\phi_{1+\alpha}\left(  \frac{s}%
{t}\right)  +(1-\alpha)\phi_{\alpha}\left(  \frac{s}{t}\right)  \right]
\label{36c}%
\end{equation}
of arguments $s,t>0$ parametrized by $\alpha\geq0.$

\paragraph{Theorem 3.1.2.}

Decomposable pseudodistances (\ref{36aa}) are for all $(P,Q)\in\mathcal{P}%
\otimes\mathcal{P}^{+}$ modifed power divergences $D_{\alpha}(P,Q)$ and
$D_{1+\alpha}(P,Q)$\ in the sense that the pseudodistance densities
$\psi_{\alpha}(p,q)$\ are\ weighted densities $\varphi_{\alpha}(p,q)$ of the
mixed power divergences
\begin{equation}
\int\varphi_{\alpha}(p,q)\,\mathrm{d}\lambda_{Q}=\alpha\,D_{1+\alpha
}(P,Q)+(1-\alpha)\,D_{\alpha}(P,Q)\label{36f}%
\end{equation}
with the power weights $w_{\alpha}(q)=q^{\alpha}$, i.e. $\psi_{\alpha
}(p,q)=w_{\alpha}(q)\varphi_{\alpha}(p,q)$ on $(\mathcal{X},\mathcal{A})$.

\paragraph{Proof.}

By (\ref{36c}),%
\begin{align}
\int\varphi_{\alpha}(p,q)\,\mathrm{d}\lambda & =\alpha\,\int\phi_{1+\alpha
}(p,q)\,\mathrm{d}\lambda+(1-\alpha)\,\int\phi_{\alpha}(p,q)\,\mathrm{d}%
\lambda\nonumber\\
& \nonumber\\
& =\alpha\,D_{1+\alpha}(P,Q)+(1-\alpha)\,D_{\alpha}(P,Q).\label{c7}%
\end{align}
By (\ref{34}), $\psi_{\alpha}(s,t)=t^{\alpha}\varphi_{\alpha}(s,t)$ so that,
by the first equality in (\ref{36aa}),%
\[
\mathfrak{D}_{\alpha}(P_{\theta},Q)=\int\psi_{\alpha}(p_{\theta}%
,q)\,\mathrm{d}\lambda=\int w_{\alpha}(q)\varphi_{\alpha}(p,q))\,\mathrm{d}%
\lambda.\medskip
\]
This together with (\ref{c7}) implies the desired result.\hfill$\blacksquare$

\bigskip

Due to Theorem 3.1.2, we call the pseudodistances $\mathfrak{D}_{\alpha}(P,Q)$
simply $\boldsymbol{power}$ $\boldsymbol{pseudo-}$ $\boldsymbol{distances}$ of
orders $\alpha\geq0$. The next theorem guarantees finiteness and continuity of
these divergences. It is restricted to the families $\mathcal{P}$\ satisfying
for some $\beta>0$ the condition%
\begin{equation}
p^{\beta},\ q^{\beta},\ln p\in\mathbb{L}_{1}(Q)\text{ \ \ for all }%
P\in\mathcal{P},\text{ }Q\in\mathcal{P}^{+}.\label{36B}%
\end{equation}

\paragraph{Theorem 3.1.3.}

If (\ref{36B}) holds for some $\beta>$ $0$\ then for all $0\leq\alpha\leq
\beta,$\ the modified power divergences are well defined by (\ref{36aa})\ and
finite, satisfying for all $P\in\mathcal{P},$ $Q\in\mathcal{P}^{+}$\ the
continuity relation%
\begin{equation}
\lim_{\alpha\downarrow0}\mathfrak{D}_{\alpha}(P,Q)=\mathfrak{D}_{0}%
(P,Q).\label{36g}%
\end{equation}

\paragraph{Proof.}

By (\ref{345}),%
\[
\mathfrak{D}_{\alpha}(P,Q)={\frac{1}{1+\alpha}}\left(  P\cdot p^{\alpha
}-Q\cdot q^{\alpha}\right)  +Q\cdot\left(  \frac{q^{\alpha}-1}{\alpha}%
-\frac{p^{\alpha}-1}{\alpha}\right)
\]
By means of the indicator function \textbf{1 }we can decompose%
\[
P\cdot p^{\alpha}=P\cdot(p^{\alpha}\text{\textbf{1}}(p\leq1))+P\cdot
(p^{\alpha}\text{\textbf{1}}(p>1))
\]
where%
\[
\lim_{\alpha\downarrow0}P\cdot(p^{\alpha}\text{\textbf{1}}(p\leq
1))=P\cdot(\text{\textbf{1}}(p\leq1))
\]
by the Lebesgue bounded convergence theorem for integrals and
\[
\lim_{\alpha\downarrow0}P\cdot(p^{\alpha}\text{\textbf{1}}(p>1))=P\cdot
(\text{\textbf{1}}(p>1))
\]
by the monotone convergence theorem for integrals. Therefore
\[
\lim_{\alpha\downarrow0}P\cdot p^{\alpha}=P\cdot(\text{\textbf{1}}%
(p\leq1))+P\cdot(\text{\textbf{1}}(p>1))=1
\]
Similarly, $\lim_{\alpha\downarrow0}Q\cdot q^{\alpha}=1$. The convergences
\[
\lim_{\alpha\downarrow0}Q\cdot\frac{q^{\alpha}-1}{\alpha}=Q\cdot\ln q\text{
\ and \ }\lim_{\alpha\downarrow0}Q\cdot\frac{p^{\alpha}-1}{\alpha}=Q\cdot\ln p
\]
follow from the monotone convergence as well, because for every fixed $t>0$\
\[
\frac{\mathrm{d}}{\mathrm{d}\alpha}\frac{t^{\alpha}-1}{\alpha}=\frac
{1-t^{\alpha}(1-\ln t)}{\alpha^{2}}\geq\frac{1-t^{\alpha}t^{-\alpha}}%
{\alpha^{2}}=0
\]
so that the expressions $(q^{\alpha}-1)/\alpha$ and $(p^{\alpha}-1)/\alpha
$\ tend monotonically to $\ln q$\ and $\ln p$.\hfill$\blacksquare$

\bigskip

By (\ref{34A}) the expressions $\mathfrak{D}_{\psi_{\alpha}}^{0}(P_{\theta}%
)$\ considered in(\ref{31c}), (\ref{31d}) are now given by%
\[
\mathfrak{D}_{\alpha}^{0}(P_{\theta})=\frac{1}{1+\alpha}\int p_{\theta
}^{1+\alpha}\,\mathrm{d}\lambda\text{ \ for all }\alpha\geq0\text{.}%
\]
Therefore the formulas (\ref{31c}), (\ref{31d}) and (\ref{28a}) lead to the
$\boldsymbol{power}$ $\boldsymbol{pseudodistance}$ $\boldsymbol{estimators}$
(briefly, $\min\mathfrak{D}_{\alpha}$-estimators)%
\begin{equation}
\theta_{\alpha,n}=\left\{
\begin{array}
[c]{ll}%
\mbox{argmin}_{\theta}\left[  \frac{{\Large 1}}{{\Large 1+\alpha}}\int
p_{\theta}^{1+\alpha}\,\mathrm{d}\lambda-{\frac{{\Large 1}}{{\Large n\alpha}}%
}\sum_{i=1}^{n}p_{{\theta}}^{\alpha}(X_{i})\right]  & \text{ \ \ \ }%
\mbox{if}\ \alpha>0\medskip\\
\mbox{argmax}_{\theta}\frac{{\Large 1}}{{\Large n}}\sum_{i=1}^{n}\ln
p_{\theta}(X_{i}) & \text{ \ \ \ }\mbox{if}\text{ }\alpha=0.
\end{array}
\right. \label{pd}%
\end{equation}
Here the upper objective function can be replaced by%
\begin{align*}
& \frac{1-\alpha}{\alpha}+\frac{1}{1+\alpha}\int p_{\theta}^{1+\alpha
}\,\mathrm{d}\lambda-{\frac{1}{n\alpha}}\sum_{i=1}^{n}p_{{\theta}}^{\alpha
}(X_{i})\\
& =\frac{1}{1+\alpha}\int p_{\theta}^{1+\alpha}\,\mathrm{d}\lambda-{\frac
{1}{n}}\sum_{i=1}^{n}\frac{p_{{\theta}}^{\alpha}(X_{i})-1}{\alpha}-1
\end{align*}
which tends for $\alpha\downarrow0$ to the lower criterion function.
Therefore, if for a fixed $n$\ the minima of all functions in (\ref{pd}) are
in a compact subset of $\Theta$\ and the MLE ${\theta}_{n,0}$\ is unique then
\begin{equation}
\lim_{\alpha\downarrow0}{\theta}_{n,\alpha}={\theta}_{n,0}.\label{pd+}%
\end{equation}
$\medskip$

\paragraph{Example 3.1.3: $L_{2}$-estimator revisited.}

By (\ref{pd}), the $\min\mathfrak{D}_{\alpha}$-estimator of order $\alpha
=1$\ is defined by%
\[
{\theta}_{1,n}=\mbox{argmin}_{{\theta}}\left[  \int p_{{\theta}}%
^{2}\,\mathrm{d}\lambda-{\frac{2}{n}}\sum_{i=1}^{n}p_{{\theta}}(X_{i})\right]
\]
so that it is nothing but the $L_{2}$-estimator ${\theta}_{n}$\ from Example
3.1.2. The family of estimators ${\theta}_{n,\alpha}$\ from (\ref{pd})
smoothly connects this robust estimator with the efficient MLE ${\theta}%
_{n,0}$ when the parameter $\alpha$\ decreases from $1$\ to $0$.$\medskip$

\paragraph{Remark 3.1.1.}

The special class of the $\min\mathfrak{D}_{\alpha}$-estimators $\theta
_{\alpha,n}$\ given by (\ref{pd}) was proposed by Basu et al. (1998) who
confirmed their efficiency for $\alpha\approx0$ and their intuitively expected
robustness for $\alpha>0$. These authors called $\theta_{\alpha,n}$
\textsl{minimum density power divergence estimators} without actual
clarification of the relation of the \textquotedblleft density power
divergences" $\mathfrak{D}_{\alpha}(P,Q)$\ to the standard power divergences
$D_{\alpha}(P,Q)$ studied in Liese and Vajda (1987) and Read and Cressie
(1988). Theorem 3.1.2 which explains $\mathfrak{D}_{\alpha}(P,Q) $\ as a
convex mixture of modified power divergences $D_{\alpha}(P,Q) $\ and
$D_{1+\alpha}(P,Q)$\ where the modification means weighting of the power
divergence densities by the power $q^{\alpha}$\ of the second probability
density, is in this respect an interesting new result.$\medskip$

\paragraph{Remark 3.1.2.}

The formula (\ref{pd}) can be given the equivalent form%
\begin{equation}
{\theta}_{\alpha,n}=\mbox{argmax}_{{\theta}}\left\{
\begin{array}
[c]{ll}%
{\frac{{\Large 1}}{{\Large n}}}\sum_{i=1}^{n}{\frac{{\Large 1}}{{\Large \alpha
}}}\left(  p_{{\theta}}^{\alpha}(X_{i})-1\right)  -\frac{1}{1+\alpha}\int
p_{{\theta}}^{1+\alpha}\,\mathrm{d}\lambda & \text{ \ \ \ }\mbox{if}\ \alpha
>0\medskip\\
{\frac{1}{n}}\sum_{i=1}^{n}\ln p_{{\theta}}(X_{i})-1 & \text{ \ \ \ }%
\mbox{if}\text{ }\alpha=0.
\end{array}
\right. \label{37}%
\end{equation}
If the integral does not depend on ${\theta}$\ then (\ref{37}) is equivalent
to%
\begin{equation}
{\theta}_{\alpha,n}=\mbox{argmax}_{{\theta}}\left\{
\begin{array}
[c]{ll}%
{\frac{{\Large 1}}{{\Large n}}}\sum_{i=1}^{n}\frac{{\Large 1}}{{\Large \alpha
}}\left(  p_{{\theta}}^{\alpha}(X_{i})-1\right)  & \text{ \ \ \ }%
\mbox{if}\ \alpha>0\medskip\\
{\frac{{\Large 1}}{{\Large n}}}\sum_{i=1}^{n}\ln p_{{\theta}}(X_{i}) & \text{
\ \ \ }\mbox{if}\text{ }\alpha=0.
\end{array}
\right. \label{37a}%
\end{equation}
This subclass of general $\min\mathfrak{D}_{\alpha}$-estimators (\ref{37}) was
included in a wider family of generalized MLE's introduced and studied
previously in Vajda (1984,1986). However, the whole class (\ref{37}) was not
introduced there.$\medskip$

If the statistical model $\langle({\mathcal{X}},{\mathcal{A}});{\mathcal{P}%
}=(P_{\theta}:\theta\in\Theta)\rangle$ is reparametrized by $\vartheta
=\vartheta(\theta)$ then the new $\min{\mathcal{D}}_{\alpha}$-estimates
$\vartheta_{\alpha_{n}}$ are related to the original $\theta_{\alpha,n}$ by
$\vartheta_{\alpha,n}=\vartheta(\theta_{\alpha,n})$. If the observations
$x\in{\mathcal{X}}$ are replaced by $y=T(x)$ where $T:({\mathcal{X}%
},{\mathcal{A}})\mapsto({\mathcal{Y}},{\mathcal{B}})$ is a measurable
statistic with the inverse $T^{-1}$ then the densities
\[
\tilde{p}_{\theta}={\frac{\mathrm{d}\tilde{P}_{\theta}}{\mathrm{\ d}%
\tilde{\lambda}}}%
\]
in the transformed model $\langle({\mathcal{Y}},{\mathcal{B}});\tilde
{{\mathcal{P}}}=(\tilde{P}_{\theta}=P_{\theta}T^{-1}:\theta\in\Theta)\rangle$
w.r.t. $\sigma$-finite dominating measure $\tilde{\lambda}=\lambda T^{-1}$ is
related to the original densities $p_{\theta}$ by
\begin{equation}
\tilde{p}_{\theta}(y)=p_{\theta}(T^{-1}y)\,{\mathcal{J}}_{T}(y)\label{c8}%
\end{equation}
where ${\mathcal{J}}_{T}(y)=\mathrm{d}\lambda T^{-1}/\mathrm{d}\tilde{\lambda
}$ is a generalized Jacobian of the statistic $T$. If ${\mathcal{X}}$,
${\mathcal{Y}}$ are Euclidean spaces, $\lambda$\ is the Lebesque measure and
the inverse mapping $\ H=T^{-1}\ $is differentiable then ${\mathcal{J}}%
_{T}(y)$ is the determinant%
\[
{\mathcal{J}}_{T}(y)=\left\vert \frac{\mathrm{d}}{\mathrm{d}y}H(y)\right\vert
.
\]

The $\min{\mathcal{D}}_{\alpha}$-estimators are in general not equivariant
w.r.t. invertible transformations of observations $T$, unless $\alpha=0$. The
following theorem generalizes similar result of Section~3.4 in Basu et al. (1998).

\paragraph{Theorem 3.1.4.}

The $\min\mathfrak{D}_{\alpha}$-estimates $\tilde{\theta}_{\alpha,n}$ in the
above considered transformed model coincide with the original $\min
\mathfrak{D}_{\alpha}$-estimates $\theta_{\alpha,n}$ if the Jacobian
${\mathcal{J}}_{T}$ of transformation is a nonzero constant on the transformed
observation space ${\mathcal{Y}}$. Thus if ${\mathcal{X}},{\mathcal{Y}}$ are
Euclidean spaces then the $\min\mathfrak{D}_{\alpha}$-estimators are
equivariant under linear statistics $Tx=ax+b$.

\paragraph{Proof.}

For $\alpha=0$ the $\min\mathfrak{D}_{\alpha}$-estimator is the MLE whose
equivariance is well known. For $\alpha>0$, by definition (\ref{pd}) and
(\ref{c8}),
\begin{align*}
\tilde{\theta}_{\alpha,n}  & =\mbox{argmin}_{\theta}\left[  {\frac{1}%
{1+\alpha}}\int_{\mathcal{Y}}\tilde{p}_{\theta}^{1+\alpha}\,\mathrm{d}%
\tilde{\lambda}-{\frac{1}{n\alpha}}\sum_{i=1}^{n}\tilde{p}_{\theta}^{\alpha
}(TX_{i})\right]  \medskip\\
& =\mbox{argmin}_{\theta}\left[  {\frac{1}{1+\alpha}}\int p_{\theta}%
^{1+\alpha}\,{\mathcal{J}}_{T}\,\mathrm{d}\lambda-{\frac{1}{n\alpha}}%
\sum_{i=1}^{n}p_{\theta}(X_{i})\,{\mathcal{J}}_{T}(TX_{i})\right]  .
\end{align*}
We see by comparison with (\ref{pd}) that $\tilde{\theta}_{\alpha,n}%
=\theta_{\alpha,n}$ if ${\mathcal{J}}_{T}$ is a nonzero constant on
${\mathcal{Y}}$. If $\alpha=0$ then the estimator is MLE and its equivariance
is well known.\hfill$\blacksquare$

\bigskip

%{\LARGE Vsuvka 2 zde}

Next we derive the influence function of the $\min\mathfrak{D}_{\alpha}%
$-estimators $\theta_{\alpha,n}$ of (\ref{pd}). Similarly as in (\ref{p1}), we
use
\[
s_{\theta}={\frac{\mathrm{d}}{\mathrm{d}\theta}}\ln p_{\theta}\quad
\mbox{and}\quad\mathring{s}_{\theta}=\left(  {\frac{\mathrm{d}}{\mathrm{\ d}%
\theta}}\right)  ^{\text{t}}s_{\theta}.
\]
It holds $\theta_{\alpha,n}=T_{\alpha}(P_{n})$ where $T_{\alpha}(Q)$ for
$Q\in{{\mathcal{Q}}}$ solves the equation $\Psi_{\alpha}(Q,\theta)\equiv
Q\cdot{\mbox{\boldmath$\psi$}}(x,\theta)=0$ for
\begin{align}
\mbox{\boldmath$\psi$}_{\alpha}(x,\theta)  & ={\frac{\mathrm{d}}%
{\mathrm{d}\theta}}\left(  {\frac{p_{\theta}^{\alpha}}{\alpha}}-{\frac
{1}{1+\alpha}}\int p_{\theta}^{1+\alpha}\,\mathrm{d}\lambda\right)
\medskip\nonumber\\
& =p_{\theta}^{\alpha}(x)\,s_{\theta}(x)-P_{\theta}\cdot p_{\theta}^{\alpha
}s_{\theta}.\label{ps}%
\end{align}
Since
\begin{equation}
\mathring{\psi}_{\alpha}(x,\theta)=\left(  {\frac{\mathrm{d}}{\mathrm{d}%
\theta}}\right)  ^{\text{t}}\,\mbox{\boldmath$\psi$}_{\alpha}(x,\theta
)=\Pi_{\alpha,\theta}(x)-P_{\theta}\cdot\left(  \Pi_{\alpha,\theta}+p_{\theta
}^{\alpha}s_{\theta}s_{\theta}^{\text{t}}\right) \label{f9}%
\end{equation}
for%
\begin{equation}
\Pi_{\alpha,\theta}=p_{\theta}^{\alpha}\,\left(  \alpha s_{\theta}s_{\theta
}^{\text{t}}+\mathring{s}_{\theta}\right)  ,\label{f10}%
\end{equation}
the matrix (\ref{f6}) is given for all $Q\in{\mathcal{P}}^{+}$ by the formula%
\begin{equation}
\mbox{\boldmath$I$}_{\alpha}(Q)=Q\cdot\Pi_{\alpha,\tau_{\alpha}}%
(x)-P_{\tau_{\alpha}}\cdot\left(  \Pi_{\alpha,\tau_{\alpha}}+p_{\tau_{\alpha}%
}^{\alpha}s_{\tau_{\alpha}}s_{\tau_{\alpha}}^{\text{t}}\right)  \text{ \ for
}\tau_{\alpha}=T_{\alpha}(Q)\in\Theta\label{f11}%
\end{equation}
In particular,
\begin{equation}
\mbox{\boldmath$I$}_{\alpha}(\theta)\equiv\mbox{\boldmath$I$}_{\alpha
}(P_{\theta})=-P_{\theta}\cdot p_{\theta}^{\alpha}s_{\theta}s_{\theta
}^{\text{t}}.\label{f12}%
\end{equation}

\bigskip

By combining (\ref{ps}), (\ref{f11}) and (\ref{f12}) with Theorem~1.1 and
Corollary 1.1, and taking into account the Fisher consistency in Theorem 3.1,
we obtain the following extension of the influence function obtained in
\S \thinspace3.3 of Basu et al. (1998) to arbitrary observation spaces
$({\mathcal{X}},{\mathcal{A}})$.

\paragraph{Theorem 3.1.5.}

If the influence function (\ref{f3}) at $Q\in{\mathcal{P}}^{+}$\ or
$P_{\theta}\in{\mathcal{P}}$\ exists for some $\min\mathfrak{D}_{\alpha}%
$-estimator $\theta_{\alpha,n}=T_{\alpha}(P_{n})$ then it is given by the
formula
\begin{equation}
\mbox{IF}(x;T_{\alpha},Q)=-\mbox{\boldmath$I$}_{\alpha}(Q)^{-1}\left[
p_{\tau_{\alpha}}^{\alpha}(x)\,s_{\tau_{\alpha}}(x)-P_{\tau_{\alpha}}\cdot
p_{\tau_{\alpha}}^{\alpha}s_{\tau_{\alpha}}\right]  \text{ \ for }\tau
_{\alpha}=T_{\alpha}(Q)\label{f13}%
\end{equation}
or%
\begin{equation}
\mbox{IF}(x;T_{\alpha},\theta)=-\mbox{\boldmath$I$}_{\alpha}(\theta
)^{-1}\left[  p_{\theta}^{\alpha}(x)\,s_{\theta}(x)-P_{\theta}\cdot p_{\theta
}^{\alpha}s_{\theta}\right] \label{f14}%
\end{equation}
respectively.

%%%%%%%%%%%%%%%%%%%%%%%%%%

\subsection{Applications in the normal family}

Consider the general normal family of Example 2.1.1. By (\ref{37}),
$\min\mathfrak{D}_{\alpha}$-estimator ${\theta}_{\alpha,n}=({\mu}_{\alpha
,n},{\sigma}_{\alpha,n})$ is the MLE given by (\ref{mle}) when $\alpha=0$.
Since
\begin{equation}
\int p_{{\theta}}^{1+\alpha}\,\mathrm{d}x=\int\left(  \frac{\exp\{-(x-{\mu
})^{2}/2{\sigma}^{2}\}}{(2\pi{\sigma}^{2})^{1/2}}\right)  ^{1+\alpha
}\,\mathrm{d}x=\frac{(1+\alpha)^{-1/2}}{(2\pi{\sigma}^{2})^{\alpha/2}%
},\label{37d}%
\end{equation}
we see from (\ref{37}) that the $\min\mathfrak{D}_{\alpha}$-estimates are for
$\alpha>0$\ given by%
\begin{align}
({\mu}_{\alpha,n},{\sigma}_{\alpha,n})  & =\mbox{argmax}_{{\mu},{\sigma}%
}\left[  {\frac{1}{\alpha n}}\sum_{i=1}^{n}\frac{\exp\left\{  -\alpha
(X_{i}-{\mu})^{2}/2{\sigma}^{2}\right\}  }{(2\pi{\sigma}^{2})^{\alpha/2}%
}-\frac{(1+\alpha)^{-3/2}}{(2\pi{\sigma}^{2})^{\alpha/2}}\right]
\medskip\medskip\nonumber\\
& =\mbox{argmax}_{{\mu},{\sigma}}\frac{1}{n\sigma^{\alpha}}\sum_{i=1}%
^{n}\left(  \exp\left\{  -\alpha\frac{(X_{i}-{\mu})^{2}}{2{\sigma}^{2}%
}\right\}  -\frac{\alpha}{(1+\alpha)^{3/2}}\right)  .\label{37e}%
\end{align}
Notice that in practical applications, the trivial "solutions" $({\mu}%
_{\alpha,n},{\sigma}_{\alpha,n})=(\max_{i}X_{i},0)$ can be avoided by
restricting the maximization to the scales bouded avay from zero.

\paragraph{Example 3.2.1: Power pseudodistance estimators of location.}

Consider the normal family $\mathcal{P}=\left\{  P_{\mu}:\mu\in\mathbb{R}%
\right\}  $\ of Example 2.1.2\ where $P_{\mu}$\ are given by the densities
$p_{\mu}(x)=p(x-\mu)$\ for the standard normal density $p(x)$. This family
satisfies the condition of the formula (\ref{37a}) so that from (\ref{pd}) or
(\ref{37a}) we obtain the $\min\mathfrak{D}_{\alpha}$-estimators ${\mu
}_{\alpha,n}=T_{\alpha}(P_{n})$ of location $\mu_{0}\in\mathbb{R}$\ in this
family given by
\begin{equation}
{\mu}_{\alpha,n}=\mbox{argmax}_{{\mu}}\left\{
\begin{array}
[c]{ll}%
\sum_{i=1}^{n}\exp\{-\alpha(X_{i}-{\mu})^{2}/2\} & \text{ \ \ \ }%
\mbox{if}\ \alpha>0\medskip\\
{-}\sum_{i=1}^{n}(X_{i}-{\mu})^{2} & \text{ \ \ \ }\mbox{if}\text{ }\alpha=0.
\end{array}
\right. \label{37b}%
\end{equation}
Equivalently, they can be obtained by inserting $\sigma=1$ in (\ref{37e}%
)$.$\ If $\alpha=0$ then ${\mu}_{\alpha,n}$ is the standard sample mean.
$\medskip$

The estimators of location (\ref{37b}) were introduced and studied as part of
larger class of estimators by Vajda (1986, 1989a,b). He proved that if the
observations are generated by $Q_{\mu_{0}}\in\mathcal{P}^{+}$ with density
$q(x-\mu_{0})$\ for unimodal $q(x)$\ symmetric about $x=0$ then these
estimators consistently estimate $\mu_{0}.$\ For $q$\ differentiable with
derivative $q^{\prime}$\ he found the influence functions
\begin{equation}
\text{IF}(x;T_{\alpha},q)=\frac{x\exp\{-\alpha x^{2}/2\}}{\int x\exp\{-\alpha
x^{2}/2\}\,q^{\prime}(x)\,\mathrm{d}x}\text{ \ \ for }\alpha\geq0.\label{37c}%
\end{equation}
This formula follows also from (\ref{f12}) and (\ref{f13}) where in this case%
\begin{equation}
s_{\mu}(x)=x-\mu\Pi_{\alpha,\mu}=p_{\mu}^{\alpha}\,\left[  \alpha\left(
x-\mu\right)  ^{2}-1\right]  \text{ \ and \ }P_{\mu}\cdot p_{\mu}^{\alpha
}s_{\mu}=0.\label{37cc}%
\end{equation}
Indeed, (\ref{37cc}) implies $P_{\mu}\cdot p_{\mu}^{\alpha}s_{\mu}=0$ and
$p_{0}^{\alpha}(x)s_{0}(x)=x\exp\{-\alpha x^{2}/2\}.\left(  2\pi\right)
^{-\alpha/2}$ so that the numerator in (\ref{37c}) follows from (\ref{f13}).
Using the identities%
\[
P_{\mu}\cdot\left(  \Pi_{\alpha,\mu}+p_{\mu}^{\alpha}s_{\mu}^{2}\right)  =\int
p_{\mu}^{1+\alpha}\,\left[  \left(  1+\alpha\right)  \left(  x-\mu\right)
^{2}-1\right]  \,\mathrm{d}x=0
\]
and
\[
\int x\,p_{0}(x)\,q^{\prime}(x)\mathrm{d}x+\int\left[  p_{0}(x)+xp_{0}%
^{\prime}(x)\right]  \,q(x)\,\mathrm{d}x=0
\]
we get from (\ref{37cc}) and (\ref{f11})%
\[
\mbox{\boldmath$I$}_{\alpha}(q)=\left(  2\pi\right)  ^{-\alpha/2}\int
x\exp\{-\alpha x^{2}/2\}\,q^{\prime}(x)\,\mathrm{d}x\,
\]
so that the denominator in (\ref{37c}) follows from (\ref{f13}).$\medskip$

The particular influence curve obtained in (\ref{37c}) for $\alpha=1/5$\ very
closely and smoothly approximates the trapezoidal IF$(x;25A,q)$ of the
estimator referred as the best under the name \textbf{Hampel's choice 25A} in
the \textsl{Princeton Robustness Study} of Andrews et al. (1972). This study
as well as the estimator of location $25A$ were influential and frequently
cited in the first decades of robust statistics. The asymptotic normality%
\[
\sqrt{n}({\mu}_{\alpha,n}-\mu_{0})\longrightarrow N(0,\sigma_{\alpha}%
^{2})\text{ \ for \ }\sigma_{\alpha}^{2}=\int\text{IF}^{2}(x;T_{\alpha
},q)q(x)\mathrm{d}x
\]
in the data generating model $Q_{\mu_{0}}$\ was established in Vajda (1986,
1989a,b) too, and the simulations presented there demonstrated that the
estimator ${T}_{1/5}$\ overperformed the set of 6 robust estimators of
location including those considered as the most prominent at that
time.$\medskip$

\paragraph{Example 3.2.2: Power pseudodistance estimators of scale.}

Consider the normal family $\mathcal{P}=\left\{  P_{\sigma}:\sigma>0\right\}
$\ of Example 2.1.3\ where $P_{\sigma}$\ are given by the densities
$p_{\sigma}(x)=p(x/\sigma)/\sigma$\ for the standard normal density $p(x)$. If
$\alpha=0$ then, by (\ref{37}), the $\min\mathfrak{D}_{\alpha}$-estimator
${\sigma}_{\alpha,n}=T_{\alpha}(P_{n})$ is the standard MLE of scale given in
(\ref{mle}). Otherwise we get from (\ref{37e}) by inserting $\mu=0$%
\begin{equation}
{\sigma}_{\alpha,n}=\mbox{argmax}_{{\sigma}}{\frac{1}{\sigma^{\alpha}\text{
}n}}\sum_{i=1}^{n}\left[  \exp\left\{  -\frac{\alpha X_{i}{}^{2}}{2{\sigma
}^{2}}\right\}  -\frac{\alpha}{(1+\alpha)^{3/2}}\right]  ,\text{ \ \ }%
\alpha>0.\label{cc}%
\end{equation}
Taking into account here
\[
\frac{1}{n}\sum_{i=1}^{n}\exp\left\{  -\frac{\alpha X_{i}{}^{2}}{2{\sigma}%
^{2}}\right\}  =\int\exp\left\{  -\frac{\alpha x{}^{2}}{2{\sigma}^{2}%
}\right\}  \mathrm{d}P_{n}(x)
\]
we find more general formula%
\[
T_{\alpha}(Q)=\text{\textrm{argmin}}_{\sigma}M_{\alpha}(Q,\sigma)\text{
\ \ for \ }Q\in\mathcal{P}^{+}%
\]
where%
\[
M_{\alpha}(Q,\sigma)={\frac{1}{\sigma^{\alpha}}}\int\exp\left\{  -\frac{\alpha
x{}^{2}}{2{\sigma}^{2}}\right\}  \mathrm{d}Q(x)-\frac{\alpha}{(1+\alpha
)^{3/2}}.
\]
By (\ref{bigpsi}) and (\ref{f00}),%
\begin{align}
\psi_{\alpha}(x,\sigma)  & =\frac{\mathrm{d}}{\mathrm{d}\sigma}M_{\alpha
}(\delta_{x},\sigma)=\frac{\mathrm{d}}{\mathrm{d}\sigma}{\frac{1}%
{\sigma^{\alpha}}}\left[  \exp\left\{  -\frac{\alpha x{}^{2}}{2{\sigma}^{2}%
}\right\}  -\frac{\alpha}{(1+\alpha)^{3/2}}\right]  \medskip\medskip
\nonumber\\
& =\frac{1}{\sigma^{1+\alpha}}\left[  \left(  \frac{x^{2}}{\sigma^{2}%
}-1\right)  \exp\left\{  -\frac{\alpha x^{2}}{2\sigma^{2}}\right\}
+\frac{\alpha}{\left(  1+\alpha\right)  ^{3/2}}\right]  .\label{cc1}%
\end{align}
The last formula will be used to evaluate the influence function. Before doing
so we shall verify it by checking the Fisher consistency condition%
\begin{equation}
P_{\sigma_{0}}\cdot\psi_{\alpha}(x,\sigma)=0\ \ \ \text{if\ and only if
\ \ }\sigma=\sigma_{0}\label{cc5}%
\end{equation}
guaranteed by Theorem 3.1. We shall use the substitutions
\begin{equation}
\sigma_{\alpha}=\frac{\sigma}{\sqrt{\alpha}},\text{ \ \ }s_{\alpha}%
=\frac{\sigma_{\alpha}\sigma_{0}}{\sqrt{\sigma_{\alpha}^{2}+\sigma_{0}^{2}}%
}\label{ccc5}%
\end{equation}
and the formula
\begin{equation}
\exp\left\{  -\frac{\alpha x^{2}}{2\sigma^{2}}\right\}  p_{\sigma_{0}}%
=\frac{s_{\alpha}}{\sigma_{0}}p_{s_{\alpha}}.\label{cccc5}%
\end{equation}
Then%
\begin{align*}
\int\left(  \frac{x^{2}}{\sigma^{2}}-1\right)  \exp\left\{  -\frac{\alpha
x^{2}}{2\sigma^{2}}\right\}  p_{\sigma_{0}}(x)\mathrm{d}x  & =\frac{s_{\alpha
}}{\sigma_{0}}\int\left(  \frac{x^{2}}{\sigma^{2}}-1\right)  p_{s_{\alpha}%
}(x)\mathrm{d}x\\
& \\
& =\frac{s_{\alpha}}{\sigma_{0}}\left(  \frac{s_{\alpha}^{2}}{\sigma^{2}%
}-1\right)  =\frac{(1-\alpha)-(\sigma/\sigma_{0})}{(\sigma/\sigma
_{0})(1+(\sigma_{0}/\sigma)\alpha)^{3/2}}%
\end{align*}
where
\[
\frac{(1-\alpha)-(\sigma/\sigma_{0})}{(\sigma/\sigma_{0})(1+(\sigma_{0}%
/\sigma)\alpha)^{3/2}}+\frac{\alpha}{\left(  1+\alpha\right)  ^{3/2}}=0
\]
if and only if $\sigma_{0}=\sigma,$ which positively verifies (\ref{cc1}%
).$\medskip$

From (\ref{cc1}) we get%
\begin{align}
\mathring{\psi}_{\alpha}(x,\sigma)  & =\frac{\mathrm{d}}{\mathrm{d}\sigma
}\mbox{\boldmath$\psi$}_{\alpha}(x,\sigma)={\frac{1}{\left(  2\pi\right)
^{\alpha/2}}.}\frac{1}{\sigma^{2+\alpha}}\exp\left\{  -\frac{\alpha x^{2}%
}{2\sigma^{2}}\right\}  \times\nonumber\\
& \nonumber\\
& \left[  \left\{  \alpha\left(  \frac{x^{4}}{\sigma^{4}}\right)
-(3+2\alpha)\left(  \frac{x^{2}}{\sigma^{2}}\right)  +1+\alpha\right\}
-\frac{\alpha}{\sqrt{1+\alpha}}\right]  .\label{cc2}%
\end{align}
Denoting for brevity as before%
\[
\tau_{\alpha}=T_{\alpha}(Q)\text{ \ \ for \ }Q\in\mathcal{P}^{+}%
\]
we obtain from (\ref{cc1}), (\ref{cc2}) and Theorem 1.1 the influence
functions of the $\min\mathfrak{D}_{\alpha}$-estimators ${\sigma}_{\alpha
,n}=T_{\alpha}(P_{n})$\ at $Q$\ for all $\alpha>0$ in the form%
\begin{align}
\text{IF}(x;T_{\alpha},Q)  & =-\frac{\psi_{\alpha}(x,\tau_{\alpha})}%
{\int\mathring{\psi}_{\alpha}(x,\tau_{\alpha})\mathrm{d}Q}\nonumber\\
& =-\frac{\sigma}{\Upsilon_{\alpha}(Q)}\left[  \exp\left\{  -\frac{\alpha
x^{2}}{2\tau_{\alpha}^{2}}\right\}  \left(  \frac{x^{2}}{\tau_{\alpha}^{2}%
}-1\right)  +\frac{\alpha}{\left(  1+\alpha\right)  ^{3/2}}\right] \label{IF}%
\end{align}
where $\Upsilon_{\alpha}(Q)$\ denotes the integral%
\begin{equation}
\int\left[  \exp\left\{  -\frac{\alpha x^{2}}{2\tau_{\alpha}^{2}}\right\}
\left\{  \alpha\left(  \frac{x}{\tau_{\alpha}}\right)  ^{4}-(3+2\alpha)\left(
\frac{x}{\tau_{\alpha}}\right)  ^{2}+1+\alpha\right\}  -\frac{\alpha}%
{\sqrt{1+\alpha}}\right]  \mathrm{d}Q.\label{cc4}%
\end{equation}
For $Q=P_{\sigma}$ the Fisher consistency implies $\tau_{\alpha}:=T_{\alpha
}(P_{\sigma})=\sigma$\ so that (\ref{IF}) and (\ref{cc4}) imply%
\[
\text{IF}(x;T_{\alpha},\sigma)=-\frac{\sigma}{\Upsilon_{\alpha}(P_{\sigma}%
)}\left[  \left(  \frac{x^{2}}{\sigma^{2}}-1\right)  \exp\left\{
-\frac{\alpha x^{2}}{2\sigma^{2}}\right\}  +\frac{\alpha}{\left(
1+\alpha\right)  ^{3/2}}\right]
\]
where the integral $\Upsilon_{\alpha}(P_{\sigma})$ reduces to\
\begin{align*}
& \int\left[  \exp\left\{  -\frac{\alpha x^{2}}{2\sigma^{2}}\right\}  \left\{
\alpha\left(  \frac{x}{\sigma}\right)  ^{4}-(3+2\alpha)\left(  \frac{x}%
{\sigma}\right)  ^{2}+1+\alpha\right\}  -\frac{\alpha}{\sqrt{1+\alpha}%
}\right]  p_{\sigma}(x)\mathrm{d}x\\
& =\frac{1}{\sqrt{1+\alpha}}\left[  \frac{3\alpha}{\left(  1+\alpha\right)
^{2}}-\frac{3+2\alpha}{1+\alpha}+1+\alpha\right]  -\frac{\alpha}%
{\sqrt{1+\alpha}}\\
& =\frac{1}{\sqrt{1+\alpha}}\left[  \frac{3\alpha-\left(  3+2\alpha\right)
\left(  1+\alpha\right)  +\left(  1+\alpha\right)  ^{2}}{\left(
1+\alpha\right)  ^{2}}\right] \\
& =-\frac{1}{\sqrt{1+\alpha}}\frac{\alpha^{2}+2}{\left(  1+\alpha\right)
^{2}}=-\frac{\alpha^{2}+2}{\left(  1+\alpha\right)  ^{5/2}}.
\end{align*}
Hence for all $\sigma>0$
\begin{equation}
\text{IF}(x;T_{\alpha},\sigma)=\frac{\left(  1+\alpha\right)  ^{5/2}\sigma
}{\alpha^{2}+2}\left[  \left(  \left(  \frac{x}{\sigma}\right)  ^{2}-1\right)
\exp\left\{  -\frac{\alpha x^{2}}{2\sigma^{2}}\right\}  +\frac{\alpha}{\left(
1+\alpha\right)  ^{3/2}}\right]  .\label{cc6}%
\end{equation}

\paragraph{Conclusion 3.2.1}

The $\min\mathfrak{D}_{\alpha}$-estimators ${\sigma}_{\alpha,n}=T_{\alpha
}(P_{n})$\ of normal scale are for all $\alpha>0$ robust in the sense that
their absolute sensitivity to the observations $x\in\mathbb{R}$\ represented
by
\[
\sup_{x\in\mathbb{R}}|\text{IF}(x;T_{\alpha},\sigma)|=\max\left\{  -\text{
IF}(0;T_{\alpha},\sigma),\text{ IF}\left(  \sigma_{\alpha};T_{\alpha}%
,\sigma\right)  \right\}  \text{ \ for }\sigma_{\alpha}=\sigma\sqrt
{\frac{2+\alpha}{\alpha}}%
\]
is bounded (cf. Hampel et al. (1986)). However, they are not insensitive
against extreme outliers because
\begin{equation}
\lim_{|x|\rightarrow\infty}\text{IF}(x;T_{\alpha},\sigma)=\text{IF}%
(\sigma;T_{\alpha},\sigma)=\frac{\alpha(1+\alpha)\sigma}{\alpha^{2}%
+2}.\label{cc7}%
\end{equation}

%%%%%%%%%%%%%%%%%%%%%%%%%%%%%%%%%%%%%%%%%%%%%%%%%

\subsection{R\'{e}nyi pseudodistance estimators}

In this subsection we propose for probability measures $P\in\mathcal{P}$ and
$Q\in\mathcal{P}^{+}$\ considered in the previous sections a family of
pseudodistances $\mathfrak{R}_{\alpha}(P,Q)$ of a R\'{e}nyi type of orders
$\alpha\geq0$ which are not of the integral type as $\mathfrak{D}_{\psi}(P,Q)$
of (\ref{25}) or $\mathfrak{D}_{\alpha}(P,Q)$\ of (\ref{36aa}). Our proposal
is based on the following theorem where%
\begin{equation}
\mathfrak{R}_{\alpha}^{0}(P)=\frac{1}{1+\alpha}\ln(P\cdot p^{\alpha})\text{
\ and \ }\mathfrak{R}_{\alpha}^{1}(Q)={\frac{1}{\alpha(1+\alpha)}}\ln(Q\cdot
q^{\alpha}).\label{38a}%
\end{equation}

%%%%%

\paragraph{Theorem 3.3.1.}

Let the condition (\ref{36B}) hold for some $\beta>0$. Then for all
$0<\alpha<\beta$%
\begin{equation}
\mathfrak{R}_{\alpha}(P,Q)=\frac{1}{1+\alpha}\ln\left(  P\cdot p^{\alpha
}\right)  +{\frac{1}{\alpha(1+\alpha)}}\ln(Q\cdot q^{\alpha})-{\frac{1}%
{\alpha}}\ln(Q\cdot p^{\alpha})\label{38}%
\end{equation}
is a family of pseudodistances decomposable in the sense%
\begin{equation}
\mathfrak{R}_{\alpha}(P,Q)=\mathfrak{R}_{\alpha}^{0}(P)+\mathfrak{R}_{\alpha
}^{1}(Q)-{\frac{1}{\alpha}}\ln(Q\cdot p^{\alpha})\label{38A}%
\end{equation}
for $\mathfrak{R}_{\alpha}^{0}(P),\mathfrak{R}_{\alpha}^{1}(Q)$ given by
(\ref{38a}), and satisfying the limit relation
\begin{equation}
\mathfrak{R}_{\alpha}(P,Q)\rightarrow\mathfrak{R}_{0}(P,Q):=Q\ln q-Q\ln
p\text{ \ \ for }\alpha\downarrow0.\label{38B}%
\end{equation}

\paragraph{Proof.}

Under (\ref{36B}), the expressions $\ln(Q\cdot q^{\alpha}),$\ $\ln(Q\cdot
p^{\alpha})$\ and $Q\cdot\ln p$\ appearing in (\ref{38}) are finite so that
the expressions $\mathfrak{R}_{\alpha}(P,Q)$\ are well defined by (\ref{38}).
Taking $\alpha>0$ and substituting%
\[
s=\frac{p^{\alpha}}{\left(  \int p^{\alpha a}\,\mathrm{d}\lambda\right)
^{1/b}},\text{ \ }t=\frac{q}{\left(  \int q^{b}\,\mathrm{d}\lambda\right)
^{1/b}}\text{ \ and \ }a=\frac{1+\alpha}{\alpha},\text{ \ }b=1+\alpha
\]
in the inequality (\ref{34a}), and integrating both sides by $\lambda$, we
obtain the H\"{o}lder inequality
\[
\int p^{\alpha}q\,\mathrm{d}\lambda\leq\left(  \int p^{1+\alpha}%
\,\mathrm{d}\lambda\right)  ^{\alpha/(1+\alpha)}\left(  \int q^{1+\alpha
}\,\mathrm{d}\lambda\right)  ^{1/(1+\alpha)}%
\]
with the equality iff $p^{\alpha a}=q^{b}$ $\ \lambda$-a.\thinspace s.$,$ i.e.
iff $p=q\ \ \lambda$-a.\thinspace s. Since the expression (\ref{38}) satisfies
for $\alpha>0$ the relation%
\begin{equation}
\mathfrak{R}_{\alpha}(P,Q)=\frac{1}{\alpha}\left\{  \ln\left[  \left(  \int
p^{1+\alpha}\,\mathrm{d}\lambda\right)  ^{\alpha/(1+\alpha)}\left(  \int
q^{1+\alpha}\,\mathrm{d}\lambda\right)  ^{1/(1+\alpha)}\right]  -\ln\int
p^{\alpha}q\,\mathrm{d}\lambda\right\}  ,\label{dkd}%
\end{equation}
we see that $\mathfrak{R}_{\alpha}(P,Q)$ is pseudodistance on the space
$\mathcal{P}$ $\otimes\mathcal{P}^{+}$. The decomposability in the sense of
(\ref{38A}) on this space is obvious and the limit relation%
\[
\mathfrak{R}_{0}(P,Q)=\lim_{\alpha\downarrow0}\mathfrak{R}_{\alpha}(P,Q)
\]
can be proved in a similar manner as in the proof of Theorem 3.1.3.\hfill
$\blacksquare$

\bigskip

There is some similarity between the decomposable pseudodistances
$\mathfrak{R}_{\alpha}(P,Q),$ $\alpha>0$\ of (\ref{38}) and the R\'{e}nyi
divergences%
\[
R_{\alpha}(P,Q)=\frac{1}{\alpha-1}\ln\left(  Q\cdot\left(  p/q\right)
^{\alpha}\right)  ,\alpha>0\text{ \ \ (cf. R\'{e}nyi (1961).}%
\]
Namely, rewriting the formula (\ref{dkd}) into the form
\[
\mathfrak{R}_{\alpha}(P,Q)=\frac{1}{\alpha+1}\ln\frac{Q\cdot\left(
p^{1+\alpha}/q\right)  }{Q\cdot p^{\alpha}}+\frac{1}{\alpha(\alpha+1)}\ln
\frac{Q\cdot q^{\alpha}}{Q\cdot p^{\alpha}}%
\]
and replacing the ratios of expectations by the expectations of ratios, we get
for $\alpha>0$ the relation%
\begin{equation}
\mathfrak{R}_{\alpha}(P,Q)=\frac{1}{\alpha+1}\ln(Q\cdot(p/q))+\frac{1}%
{\alpha(\alpha+1)}\ln(Q\cdot(q/p)^{\alpha})=\frac{1}{\alpha+1}R_{\alpha
+1}(Q,P)\label{39}%
\end{equation}
which can be extended to $\alpha=0$\ by taking on both sides the limits for
$\alpha\downarrow0$. Therefore the decomposable pseudodistances (\ref{38}) are
modified R\'{e}nyi divergences and as such, they are called
$\boldsymbol{R\acute{e}nyi}$ $\boldsymbol{pseudodistances}$. $\medskip$

Similarly as earlier in this section, we are interested in the estimators
obtained by replacing the hypothetical distribution $P_{\theta_{0}}$ in the
$\mathfrak{R}_{\alpha}$-pseudodistances $\mathfrak{R}_{\alpha}(P_{{\theta}%
},P_{\theta_{0}})$ by the empirical distribution $P_{n}$. In other words, we
are interested in the family of $\boldsymbol{R\acute{e}nyi}$
$\boldsymbol{pseudodistance\mathbf{\ }estimators}$ of orders $0\leq\alpha
\leq\beta$ (in symbols, $\min\mathfrak{R}_{\alpha}$-estimators)\ defined as
${\theta}_{n,\alpha}=T_{\alpha}(P_{n})$ for $T_{\alpha}(Q)\in\Theta$ with
$Q\in\mathcal{Q}=\mathcal{P}^{+}\cup\mathcal{P}_{\text{emp}}$ satisfying the
condition%
\begin{equation}
T_{\alpha}(Q)=\left\{
\begin{array}
[c]{ll}%
\arg\min_{{\theta}}\frac{1}{1+\alpha}\ln\left(  P_{\theta}\cdot p_{\theta
}^{\alpha}\right)  -{\frac{1}{\alpha}}\ln(Q\cdot p_{\theta}^{\alpha}) & \text{
\ \ \ }\mbox{if}\ 0<\alpha\leq\beta\medskip\\
\arg\min_{{\theta}}-\ln Q\cdot p_{{\theta}} & \text{ \ \ \ }\mbox{if}\text{
}\alpha=0.
\end{array}
\right. \label{d1}%
\end{equation}
The upper formula is for
\begin{equation}
C_{\theta}(\alpha)=\left(  P_{\theta}\cdot p_{\theta}^{\alpha}\right)
^{\alpha/(1+\alpha)}\equiv\left(  \int p_{\theta}^{1+\alpha}\mathrm{d}%
\lambda\right)  ^{\alpha/(1+\alpha)}\label{d2}%
\end{equation}
equivalent to%
\begin{equation}
T_{\alpha}(Q)=\arg\max_{\theta}M_{\alpha}(Q,\theta)\text{ \ \ for \ }%
M_{\alpha}(Q,\theta)=\frac{Q\cdot p_{\theta}^{\alpha}}{C_{\theta}(\alpha
)}\label{d3}%
\end{equation}
Alternatively, we can write
\begin{equation}
{\theta}_{n,\alpha}=\left\{
\begin{array}
[c]{ll}%
\arg\max_{\theta}C_{\theta}{\normalsize (\alpha)}^{-1}\frac{{\Large 1}%
}{{\Large n}}\sum_{i=1}^{n}p_{{\theta}}^{\alpha}(X_{i}) & \text{
\ \ \ }\mbox{if}\ 0<\alpha\leq\beta\medskip\\
\arg\max_{\theta}\frac{{\Large 1}}{{\Large n}}\sum_{i=1}^{n}\ln p_{{\theta}%
}(X_{i}) & \text{ \ \ \ }\mbox{if}\text{ }\alpha=0.
\end{array}
\right. \label{40}%
\end{equation}
For $\alpha\approx0\downarrow0$ the approximations $C_{\theta}%
{\normalsize (\alpha)}\approx1$ and%
\[
\frac{1}{\alpha}\left(  \frac{{\normalsize 1}}{{\normalsize n}}\sum_{i=1}%
^{n}p_{{\theta}}^{\alpha}(X_{i})-1\right)  =\frac{{\normalsize 1}%
}{{\normalsize n}}\sum_{i=1}^{n}\frac{p_{{\theta}}^{\alpha}(X_{i}%
)-1}{{\normalsize \alpha}}\approx\frac{{\normalsize 1}}{{\normalsize n}}%
\sum_{i=1}^{n}\ln p_{{\theta}}^{\alpha}(X_{i})
\]
indicate that the upper criterion function in (\ref{40}) tends to the lower
MLE criterion for $\alpha\downarrow0$. If $C_{\theta}(\alpha)$ does not depend
on ${\theta}$\ then the $\min\mathfrak{R}_{\alpha}$-estimates reduce to the
$\min$ $\mathfrak{D}_{\alpha}$-estimates considered in (\ref{37a}) of Remark
3.1.2, i.e.,%
\begin{equation}
{\theta}_{\alpha,n}=\mbox{argmax}_{{\theta}}\left\{
\begin{array}
[c]{ll}%
{\frac{{\Large 1}}{{\Large n}}}\sum_{i=1}^{n}p_{{\theta}}^{\alpha}(X_{i}) &
\text{ \ \ \ }\mbox{if}\ 0<\alpha<\beta\medskip\\
{\frac{{\Large 1}}{{\Large n}}}\sum_{i=1}^{n}\ln p_{{\theta}}(X_{i}) & \text{
\ \ \ }\mbox{if}\text{ }\alpha=0.
\end{array}
\right. \label{40a}%
\end{equation}
If the extremal points of all functions in (\ref{40}) are in a compact set of
$\Theta$\ then
\begin{equation}
\lim_{\alpha\downarrow0}{\theta}_{n,\alpha}={\theta}_{n,0}.\label{41}%
\end{equation}
$\medskip$

In the next theorem and its proof we use the auxiliary expressions
\[
s_{\theta}={\frac{\mathrm{d}}{\mathrm{d}\theta}}\ln p_{\theta},\quad
\mathring{s}_{\theta}=\left(  {\frac{\mathrm{d}}{\mathrm{d}\theta}}\right)
^{\text{t}}s_{\theta}\text{ \ \ \ \ \ \ \ \ \ \ \ \ \ \ \ (cf. (\ref{p1}))}%
\]
and
\[
c_{\theta}(\alpha)=\frac{\int p_{\theta}^{1+\alpha}s_{\theta}\,\mathrm{d}%
\lambda}{\int p_{\theta}^{1+\alpha}\,\mathrm{d}\lambda},\quad\mathring
{c}_{\theta}(\alpha)=\left(  {\frac{\mathrm{d}}{\mathrm{d}\theta}}\right)
^{\text{t}}c_{\theta}(\alpha)\quad\mbox{and}\quad\tau_{\alpha}=T_{\alpha
}(Q)\text{.}%
\]

\paragraph{Theorem 3.3.3.}

If the influence function (\ref{f3}) at $Q\in{\mathcal{P}}^{+}$\ or
$P_{\theta}\in{\mathcal{P}}$\ exists for some $\min\mathfrak{R}_{\alpha}%
$-estimator $\theta_{\alpha,n}=T_{\alpha}(P_{n})$ then it is given by the
formula
\begin{equation}
\mbox{IF}(x;T_{\alpha},Q)=-\mbox{\boldmath$I$}_{\alpha}(Q)^{-1}\left[
p_{\tau_{\alpha}}(x)\left(  s_{\tau_{\alpha}}(x)-c_{\tau_{\alpha}}%
(\alpha)\right)  \right] \label{ff13}%
\end{equation}
or%
\begin{equation}
\mbox{IF}(x;T_{\alpha},\theta)=-\mbox{\boldmath$I$}_{\alpha}(\theta
)^{-1}\left[  p_{\theta}(x)\left(  s_{\theta}(x)-c_{\theta}(\alpha)\right)
\right] \label{fff13}%
\end{equation}
for the matrices
\begin{equation}
\mbox{\boldmath$I$}_{\alpha}(Q)=\int\left[  \mathring{s}_{\tau_{\alpha}%
}-\mathring{c}_{\tau_{\alpha}}(\alpha)-\alpha p_{\tau_{\alpha}}^{\alpha
}\left(  s_{\tau_{\alpha}}-c_{\tau_{\alpha}}(\alpha)\right)  \left(
s_{\tau_{\alpha}}-c_{\tau_{\alpha}}(\alpha)\right)  ^{\text{t}}\right]
p_{\tau_{\alpha}}^{\alpha}\mathrm{d}Q\label{ff12}%
\end{equation}
or%
\begin{equation}
\mbox{\boldmath$I$}_{\alpha}(\theta)=\int\left[  \mathring{s}_{\theta
}-\mathring{c}_{\theta}(\alpha)-\alpha p_{\theta}^{\alpha}\left(  s_{\theta
}-c_{\theta}(\alpha)\right)  \left(  s_{\theta}-c_{\theta}(\alpha)\right)
^{\text{t}}\right]  p_{\theta}^{1+\alpha}\mathrm{d}\lambda\label{fff12}%
\end{equation}
respectively.

\paragraph{Proof.}

By (\ref{d3}), $T_{\alpha}(Q)$ for $Q\in{{\mathcal{Q}}}$ minimizes
$Q\cdot\left(  p_{\theta}^{\alpha}/C_{\theta}(\alpha)\right)  ,$ i.e. solves
the equation $\Psi_{\alpha}(Q,\theta)\equiv Q\cdot{\mbox{\boldmath$\psi$}}%
(x,\theta)=0$ for
\begin{equation}
\mbox{\boldmath$\psi$}_{\alpha}(x,\theta)\equiv\Psi_{\alpha}(\delta_{x}%
,\theta)={\frac{\mathrm{d}}{\mathrm{d}\theta}}\frac{p_{\theta}^{\alpha}%
}{C_{\theta}(\alpha)}=\frac{\alpha p_{\theta}^{\alpha}\left(  s_{\theta
}-c_{\theta}(\alpha)\right)  }{C_{\theta}(\alpha)}.\label{psi}%
\end{equation}
Further,%
\[
C_{\theta}(\alpha):=\left(  {\frac{\mathrm{d}}{\mathrm{d}\theta}}\right)
^{\text{t}}\,C_{\theta}(\alpha)=\alpha C_{\theta}(\alpha)c_{\theta}^{\text{t}%
}(\alpha)
\]
so that
\begin{align*}
\mathring{\psi}_{\alpha}(x,\theta)  & =\left(  {\frac{\mathrm{d}}%
{\mathrm{d}\theta}}\right)  ^{\text{t}}\,\mbox{\boldmath$\psi$}_{\alpha
}(x,\theta)\\
& \\
& =\frac{C_{\theta}(\alpha)\left[  \alpha^{2}p_{\theta}^{\alpha}s_{\theta
}^{\text{t}}\left(  s_{\theta}-c_{\theta}(\alpha)\right)  +\alpha p_{\theta
}^{\alpha}\left(  \mathring{s}_{\theta}-\mathring{c}_{\theta}(\alpha)\right)
\right]  -\alpha p_{\theta}^{\alpha}\left(  s_{\theta}-c_{\theta}%
(\alpha)\right)  C_{\theta}(\alpha)}{C_{\theta}(\alpha)}\\
& \\
& =\frac{\alpha^{2}p_{\theta}^{\alpha}s_{\theta}^{\text{t}}\left(  s_{\theta
}-c_{\theta}(\alpha)\right)  +\alpha p_{\theta}^{\alpha}\left(  \mathring
{s}_{\theta}-\mathring{c}_{\theta}(\alpha)\right)  -\alpha^{2}p_{\theta
}^{\alpha}s_{\theta}^{\text{t}}\left(  s_{\theta}-c_{\theta}(\alpha)\right)
c_{\theta}^{\text{t}}(\alpha)}{C_{\theta}(\alpha)}.\medskip
\end{align*}
Therefore the matrix (\ref{f6}) is given for all $Q\in{\mathcal{P}}^{+}$ by
the formula (\ref{ff12}) and (\ref{f66}) is given for $P_{\theta}%
\in\mathcal{P}$ by (\ref{fff12}). The rest is clear from Theorems 1.1 and 3.1,
and from Corollary 1.1.\hfill$\blacksquare$

%%%%%%%%%%%%%%%%%%%%%%%%%%

\subsection{Applications in the normal family}

Consider the general normal family of Example 2.1.1 for which the condition
(\ref{36B}) is satisfied for all $\beta>0$ and (\ref{37d}) implies
\begin{equation}
C_{\mu,\sigma}(\alpha)=C_{\sigma}(\alpha)=\left(  \frac{(1+\alpha)^{-1/2}%
}{(2\pi{\sigma}^{2})^{\alpha/2}}\right)  ^{{\large \alpha/}\left(
1+\alpha\right)  }=\frac{\sigma^{-\alpha^{2}/(1+\alpha)}}{c(\alpha)}\label{d4}%
\end{equation}
for all\ $\mu\in\mathbb{R}$ and the function%
\[
c(\alpha)=\left[  \left(  1+\alpha\right)  \left(  2\pi\right)  ^{\alpha
}\right]  ^{\alpha/2(1+\alpha)},\alpha>0.
\]
By (\ref{40}), the $\min\mathfrak{R}_{\alpha}$-estimator ${\theta}_{\alpha
,n}=({\mu}_{\alpha,n},{\sigma}_{\alpha,n})$ is the standard estimator of
location and scale given by (\ref{mle}) if $\alpha=0$. For $\alpha>0$\ we can
use\ the relation%
\[
\frac{\sigma^{\alpha^{2}/(1+\alpha)}}{\sigma^{\alpha}}=\sigma^{-\alpha
/(1+\alpha)}%
\]
to get from (\ref{40}) and (\ref{d4}) the highly nonstandard estimator%
\begin{equation}
({\mu}_{\alpha,n},{\sigma}_{\alpha,n})=\mbox{argmax}_{{\mu},{\sigma}}\left[
\frac{c_{\alpha}}{n{\sigma}^{\alpha/(1+\alpha)}}\sum_{i=1}^{n}\exp\left\{
-\alpha\frac{(X_{i}-{\mu})^{2}}{2{\sigma}^{2}}\right\}  \right] \label{42}%
\end{equation}
which in general differs from the $\min\mathfrak{D}_{\alpha}$-estimator
(\ref{37e}) as it will be seen in the submodel of scale below. Similarly as in
the case of power pseudodistance estimator (\ref{37e}),\ the trivial
"solutions" $({\mu}_{\alpha,n},{\sigma}_{\alpha,n})=(\max_{i}X_{i},0)$ can be
avoided in practical applications by restricting the maximization to the
scales bouded avay from zero.$\medskip$

The next example of the submodel of location illustrates the situation where
these two estimators coincide. Obviously, the constants $c_{\alpha}%
=c(\alpha)/(2\pi)^{\alpha/2}\ $play no role in the maximization and can be
replaced by 1.\medskip

\paragraph{Example 3.4.1: R\'{e}nyi pseudodistance estimators of location.}

The normal family of location introduced in Example 2.1.2\ satisfies the
condition of the formula (\ref{37a}) so that from (\ref{pd}) or (\ref{37a}) we
obtain the same $\min\mathfrak{R}_{\alpha}$-estimators ${\mu}_{\alpha,n}$\ of
location $\mu_{0}\in\mathbb{R}$\ as in (\ref{37b}). Thus to these estimators
applies all what was seen in Example 3.3.1.

\paragraph{Example 3.4.2: R\'{e}nyi pseudodistance estimators of scale.}

Consider the normal model of scale introduced in Example 2.1.3. If $\alpha=0 $
then, by (\ref{37}), the $\min\mathfrak{R}_{\alpha}$-estimator ${\sigma
}_{\alpha,n}=T_{\alpha}(P_{n})$ is the standard MLE of scale given in
(\ref{mle}). Otherwise by (\ref{42}),
\begin{equation}
{\sigma}_{\alpha,n}=\mbox{argmax}_{{\sigma}}\left[  \frac{c_{\alpha}}%
{n{\sigma}^{\alpha/(1+\alpha)}}\sum_{i=1}^{n}\exp\left\{  -\alpha\frac{X_{i}%
{}^{2}}{2{\sigma}^{2}}\right\}  \right]  ,\text{ \ \ }\alpha>0\text{
\ \ \ \ (cf. (\ref{42})).}\label{43}%
\end{equation}
It is easy to see e.g. by putting $n=1$ and $\alpha X^{2}=2$ that these
estimates differ from the $\mathfrak{D}_{\alpha}$-estimates of scale given
in(\ref{cc}). Here (\ref{d3}) for the Dirac $\delta_{x}$ implies%
\[
M_{\alpha}(\delta_{x},\sigma)=\frac{p_{\sigma}^{\alpha}(x)}{C_{\sigma}%
(\alpha)}=\frac{c_{\alpha}}{\sigma^{\alpha/(1+\alpha)}}\exp\left\{
-\frac{\alpha x^{2}}{2\sigma^{2}}\right\}
\]
and by (\ref{bigpsi}) and (\ref{f00}),%
\begin{align}
\mathbf{\psi}_{\alpha}(x,\sigma)  & =\frac{\mathrm{d}}{\mathrm{d}\sigma
}M_{\alpha}(\delta_{x},\sigma)=c_{\alpha}\frac{\mathrm{d}}{\mathrm{d}\sigma
}\left[  {\sigma}^{-\alpha/\left(  1+\alpha\right)  }\exp\left\{
-\frac{\alpha x^{2}}{2\sigma^{2}}\right\}  \right] \nonumber\\
& \nonumber\\
& =\frac{c_{\alpha}}{\sigma^{\alpha/\left(  1+\alpha\right)  }}\left[
\frac{\alpha x^{2}}{\sigma^{3}}-\frac{\alpha}{1+\alpha}\frac{1}{{\sigma}%
}\right]  \exp\left\{  -\frac{\alpha x^{2}}{2\sigma^{2}}\right\} \nonumber\\
& \nonumber\\
& =\frac{\alpha c_{\alpha}}{\sigma^{1+\alpha/\left(  1+\alpha\right)  }%
}\left[  \left(  \frac{x}{\sigma}\right)  ^{2}-\frac{1}{1+\alpha}\right]
\exp\left\{  -\frac{\alpha x^{2}}{2\sigma^{2}}\right\}  .\label{psir}%
\end{align}
This formula can be verified by checking the Fisher consistency known in
general from Theorem 3.1. Using the formulas (\ref{ccc5}) and (\ref{cccc5}) we
find%
\[
\int\left[  \left(  \frac{x}{\sigma}\right)  ^{2}-\frac{1}{1+\alpha}\right]
\exp\left\{  -\frac{\alpha x^{2}}{2\sigma^{2}}\right\}  p_{\sigma_{0}%
}(x)\mathrm{d}x=\frac{\sigma}{\sqrt{\sigma^{2}+\alpha\sigma_{0}^{2}}}\left[
\left(  \frac{\sigma_{0}^{2}}{\sigma^{2}+\alpha\sigma_{0}^{2}}\right)
^{2}-\frac{1}{1+\alpha}\right]  .
\]
Since the right-hand side is zero if and only if $\sigma=\sigma_{0},$ the
verification is positive.\bigskip

From (\ref{psir}) we evaluate after some effort the derivative%
\begin{align}
\mathbf{\mathring{\psi}}_{\alpha}(x,\sigma)  & =\frac{\mathrm{d}}%
{\mathrm{d}\sigma}\mathbf{\psi}_{\alpha}(x,\sigma)=\frac{\mathrm{d}%
}{\mathrm{d}\sigma}\frac{c_{\alpha}}{\sigma^{1+\alpha/\left(  1+\alpha\right)
}}\exp\left\{  -\frac{\alpha x^{2}}{2\sigma^{2}}\right\}  \left[
\alpha\left(  \frac{x}{\sigma}\right)  ^{2}-\frac{\alpha}{1+\alpha}\right]
\nonumber\\
& \nonumber\\
& =\frac{\alpha c_{\alpha}}{\sigma^{2+\alpha/\left(  1+\alpha\right)  }}%
\exp\left\{  -\frac{\alpha x^{2}}{2\sigma^{2}}\right\}  \eta_{\alpha}\left(
\frac{x}{\sigma}\right) \label{psidot}%
\end{align}
where
\[
\eta_{\alpha}\left(  \frac{x}{\sigma}\right)  =\left[  \alpha\left(  \frac
{x}{\sigma}\right)  ^{4}-\frac{5\alpha+3}{1+\alpha}\left(  \frac{x}{\sigma
}\right)  ^{2}+\frac{2\alpha+1}{\left(  1+\alpha\right)  ^{2}}\right]  .
\]
Thus, denoting for brevity%
\[
\tau_{\alpha}=T_{\alpha}(Q)\text{ \ \ for \ }Q\in\mathcal{P}^{+}%
\]
we obtain from (\ref{psir}), (\ref{psidot}) and Theorem 1.1 the influence
functions of the $\min\mathfrak{D}_{\alpha}$-estimators ${\sigma}_{\alpha
,n}=T_{\alpha}(P_{n})$\ at $Q$\ given for all $\alpha>0$ by%
\begin{align}
\text{IF}(x;T_{\alpha},Q)  & =-\frac{\psi_{\alpha}(x,\tau_{\alpha})}%
{\int\mathring{\psi}_{\alpha}(x,\tau_{\alpha})\mathrm{d}Q}\nonumber\\
& =-\frac{\alpha}{\Upsilon_{\alpha}(Q)}\left[  \left(  \left(  \frac{x}%
{\tau_{\alpha}}\right)  ^{2}-\frac{1}{1+\alpha}\right)  \exp\left\{
-\frac{\alpha x^{2}}{2\tau_{\alpha}^{2}}\right\}  \right] \label{ips}%
\end{align}
where%
\[
\Upsilon_{\alpha}(Q)=\int\eta_{\alpha}\left(  \frac{x}{\sigma}\right)
\exp\left\{  -\frac{\alpha x^{2}}{2\tau_{\alpha}^{2}}\right\}  \mathrm{d}Q.
\]
In the special case $Q=P_{\sigma}$\ the Fisher consistency implies that
$\tau_{\alpha}:=T_{\alpha}(P_{\sigma})=\sigma$. We use the relation%
\[
\exp\left\{  -\frac{\alpha x^{2}}{2\sigma^{2}}\right\}  p_{\sigma
}(x)=p_{\sigma_{\alpha}}(x)\frac{1}{\sqrt{1+\alpha}}\text{ \ \ for \ }%
\sigma_{\alpha}=\frac{\sigma}{\sqrt{1+\alpha}}%
\]
to obtain%
\begin{align*}
\Upsilon_{\alpha}(P_{\sigma})  & =\frac{1}{\sqrt{1+\alpha}}\int\eta_{\alpha
}\left(  \frac{x}{\sigma}\right)  p_{\sigma_{\alpha}}(x)\mathrm{d}x\\
& \\
& =\frac{1}{\left(  1+\alpha\right)  ^{1/2}}\left[  \alpha\left(  \frac
{\sigma_{\alpha}}{\sigma}\right)  ^{4}-\frac{5\alpha+3}{1+\alpha}\left(
\frac{\sigma_{\alpha}}{\sigma}\right)  ^{2}+\frac{2\alpha+1}{\left(
1+\alpha\right)  ^{2}}\right] \\
& \\
& =\frac{1}{\left(  1+\alpha\right)  ^{5/2}}\left[  3\alpha-\left(
5\alpha+3\right)  +2\alpha+1\right]  =-\frac{2}{\left(  1+\alpha\right)
^{5/2}}%
\end{align*}
independently of $\sigma>0.$ Therefore at the normal location $P_{\sigma}$\ we
get for all $\sigma>0$\ the influence functions%
\begin{equation}
\text{IF}(x;T_{\alpha},P_{\sigma})=\frac{\left(  1+\alpha\right)  ^{5/2}%
\sigma}{2}\left[  \left(  \left(  \frac{x}{\sigma}\right)  ^{2}-\frac
{1}{1+\alpha}\right)  \exp\left\{  -\frac{\alpha x^{2}}{2\sigma^{2}}\right\}
\right]  .\label{infl}%
\end{equation}
It is easy to verify that this is the influence function also in the MLE case
$\alpha=0$.

\paragraph{Conclusion 3.4.1.}

The $\min\mathfrak{R}_{\alpha}$-estimators ${\sigma}_{\alpha,n}=T_{\alpha
}(P_{n})$\ of normal scale are for all $\alpha>0$ robust in the sense that
their influence functions are bounded. They are more robust against distant
outliers than the corresponding $\min\mathfrak{D}_{\alpha}$-estimators studied
in the Subsections 3.1 and 3.2 because
\begin{equation}
\lim_{|x|\rightarrow\infty}\text{IF}(x;T_{\alpha},P_{\sigma})=0\text{
\ \ \ (cf. (\ref{ips})).}\label{ccccr}%
\end{equation}

\paragraph{Problem 3.4.1.}

Compare by simulations the mean squared errors of the $\min\mathfrak{D}%
_{\alpha}$-estimators and $\min\mathfrak{R}_{\alpha}$-estimators of location
in contaminated normal scale models%
\begin{equation}
(1-\varepsilon)P_{\sigma}+\varepsilon Q_{\sigma}\label{simul}%
\end{equation}
for%
\begin{equation}
0<\varepsilon<1/2\text{ \ and \ }Q\in\left\{  P_{3},P_{10},\text{ Logistic,
Cauchy}\right\}  .\label{model}%
\end{equation}
Verify in this manner the stronger robustness of the $\min\mathfrak{R}%
_{\alpha}$-estimators theoretically justified in the Conclusion 3.4.1.

%%%%%%%%%%%%%%%%%%%%%%%%%%

%%%%%%%%%%%%%%%%%%%%%%%%%%%

\paragraph{Acknowledgement}

This research was supported by the grants GA\thinspace\v{C}R 102/07/1131 and
M\v{S}MT 1M 0572. The authors thank the PhD student Iva Fr\'{y}dlov\'{a} for
careful reading and corrections of many previous versions of the first two
sections. They thank also to the MSc student Radim Demut for simulations of
the R\'{e}nyi estimators in contaminated families. The very promising results
obtained by him encouraged the theoretic research presented here.

%%%%%%%%%%%%%%%%%%%%%%%%%%%%%%%%%%%%%%%%%%%%%%%%%

\end{document}